\documentclass[11pt,a4paper]{amsart}
\usepackage[latin1]{inputenc}
\usepackage[arrow, matrix, curve]{xy}
\usepackage{amssymb}
\usepackage{amsmath}
\usepackage{latexsym}
\usepackage{xypic}
\usepackage{epic}
\usepackage{eepic}
\usepackage{graphicx}
\usepackage{longtable}
\usepackage{epsfig}
\usepackage{xypic}

\newtheorem{prop}{Proposition}[section]
\newtheorem{lemma}{Lemma}[section]
\newtheorem{cor}{Corollary}[section]
\newtheorem{defi}{Definition}[section]

\newcommand{\lra}{\longrightarrow}

\newcommand{\bprf}{{\it Proof.~}}
\newcommand{\eprf}{\hfill $\square$ \bigskip\par}
\newcommand{\PP}{ \mathbb{P}}
\newcommand{\C }{ \mathbb{C}}

\newcommand{\Z}{\mathbb{Z}}
\newcommand{\Q}{\mathbb{Q}}
\newcommand{\cl}{\mathcal{L}}

\input epsf
\input cyracc.def

\makeatletter
\def\blfootnote{\xdef\@thefnmark{}\@footnotetext}

\begin{document}


\title{Projective models of K3 surfaces with an even set}
\author{Alice Garbagnati and Alessandra Sarti}
\address{Alice Garbagnati, Dipartimento di Matematica, Universit\`a di Milano,
  via Saldini 50, I-20133 Milano, Italia}

\email{garba@mat.unimi.it}

\address{Alessandra Sarti, Institut f\"ur  Mathematik, Universit\"at Mainz,
Staudingerweg 9, 55099 Mainz, Germany.}

\email{sarti@mathematik.uni-mainz.de}

\begin{abstract}
The aim of this paper is to describe algebraic K3 surfaces with an
even set of rational curves or of nodes. Their  minimal possible
Picard number is nine. We completely classify these K3
surfaces and after a carefull analysis of the divisors contained in
the Picard lattice we study their projective models, giving
necessary and sufficient conditions to have an even set. Moreover
we investigate their relation with K3 surfaces with a Nikulin
involution.

\end{abstract}

\maketitle


\blfootnote {The second author was partially supported by DFG Research Grant
SA 1380/1-2.}  \blfootnote {{\it Key
words:} K3 surfaces, even sets of curves, moduli.}

\section{Introduction}\label{intro}
It is a classical problem in algebraic geometry to determine when a
set of $(-2)$-rational curves on a surface is even. This means the
following: let $L_1,\ldots,L_N$ be rational curves on a surface
$X$ then they form an {\it even set} if there is $\delta\in
Pic(X)$ such that
$$
L_1+\ldots+L_n\sim 2\delta.
$$
This is equivalent to the existence of a double cover of $X$ branched on $L_1+\ldots+L_n$. This problem
is related to the study of even sets of nodes, in fact a set of nodes
is even if the $(-2)$-rational curves in the minimal resolution are an even set. In particular the study of even sets
on surfaces plays an important role in determining the maximal number of nodes a surface can have (cf. e.g. \cite{beauv}, \cite{jaffe}). Here we restrict our attention to K3 surfaces.\\
In a famous paper of 1975 \cite{nik} Nikulin shows that an even
set of disjoint rational curves (resp. of distinct nodes) on a K3
surface contains $0$, $8$ or $16$ rational curves (nodes).
If the even set on the K3 surface $X$ is made up by sixteen
rational curves, the surface covering $X$ is birational to a
complex torus $A$ and $X$ is the Kummer surface of $A$. This
situation is studied by Nikulin in \cite{nik}. If the even set on
$X$ is made up by eight rational curves then the surface covering
$X$ is also a K3 surface. There are some more general results
about even sets of curves not necessarily disjoint. More recently
in \cite{barth1} Barth studies the case of even sets of rational
curves on quartic surfaces (i.e. K3 surfaces in $\mathbb{P}^3$)
also in the case that the curves meet each other, he finds sets containing six or ten lines too.\\
In the paper \cite{barth2} he discusses some particular even sets
of disjoint lines and nodes on K3 surfaces whose projective models
are a double cover of the plane, a quartic in $\PP^3$ or a double
cover of the quadric $\PP^1\times \PP^1$, and he gives necessary
and sufficient
conditions to have an even set.\\
Our purpose is to study algebraic K3 surfaces admitting an even
set of eight disjoint rational curves. We investigate their Picard
lattices, moduli spaces and projective models. The minimal
possible Picard number is nine, and we restrict our study to the
surfaces with this Picard number. The techniques used by Barth in
his article are mostly geometric, here we investigate first the
Picard lattices of the K3 surfaces and the ampleness of certain
divisors, then we study the projective models. We find again the
cases studied by Barth and we  discuss many new cases, with a
special attention to complete intersections. We give also an
explicit relation between the Picard lattice of an algebraic K3
surface with an even set and the Picard lattice of the K3 surface
which is its double cover. More precisely if $X$ admits an even
set of eight disjoint rational curves, then by \cite{nik}, it is
the desingularization of the quotient of a K3 surface by a Nikulin
involution (i.e.\ a symplectic automorphism of order two). The
Nikulin involutions are well known and are studied by Morrison in
\cite{morrison} and by van Geemen and Sarti in \cite{bertio}. In
\cite{bertio} the authors describe also some geometric properties
of the quotient
 by a Nikulin involution and so of
K3 surfaces with an even set of eight nodes. In \cite{nik} Nikulin
proves that a sufficient condition on a K3 surface to be a
Kummer surface (and so to have an even set made up by sixteen
disjoint rational curves) is that a particular lattice (the so
called Kummer lattice) is primitively embedded in the N\'eron
Severi group of the surface. Here we prove a similar result: a
sufficient condition on a K3 surface $X$ to be the
desingularization of the quotient of another K3 surface with a
Nikulin involution (and so to have an even set made up by eight
disjoint rational curves) is that a particular lattice (the so
called Nikulin lattice) is primitively embedded in the N\'eron
Severi group of $X$. This result is essential to describe the coarse moduli space of a K3 surface with an even set of eight disjoint rational curves.\\
In the Section \ref{section:evenset} we recall same known results
on even sets on surfaces, in particular on K3 surfaces. In the
Section \ref{section:evennik} we study algebraic K3 surfaces $X$
with Picard number nine. If $X$ admits an even set of eight
disjoint rational curves, then its N\'eron Severi group has rank
at least nine (it has to contain the eight rational curves of the
even set and a polarization, because the K3 surface is algebraic).
The main results of this section (and also two of the main results
of this paper) are the complete description of the possible
N\'eron Severi groups of rank nine of algebraic K3 surfaces
admitting an even set and the complete description of the coarse
moduli space of the algebraic K3 surfaces with an even set of
eight disjoint rational curves. Denote by $N$ the Nikulin lattice,
the minimal primitive sublattice of $H^2(X,\Z)$ containing the
$(-2)$-rational curves (cf. Definition \ref{latticenik}), then:

\smallskip

{\bf Proposition \ref{proposition: possible NS(X)}}.
{\it Let $X$ be an algebraic K3 surface with an even set of eight disjoint
rational curves and with Picard number nine, let $L$ be a divisor
generating  $N^{\perp}\subset NS(X)$, $L^2>0$. Let $d$ be a
positive integer with $L^2=2d$ and let
\begin{eqnarray*}
\mathcal{L}_{2d}=\Z L\oplus N.
\end{eqnarray*}
Then\\
(1) if $L^2\equiv 2$$\mod 4$ then $NS(X)=\mathcal{L}_{2d}$,\\
(2) if $L^2\equiv 0$$\mod 4$ then either $NS(X)=\mathcal{L}_{2d}$
or $NS(X)=\mathcal{L}_{2d}'$, where $\cl_{2d}'$ is generated by
$\cl_{2d}$ and by a class $(L/2,v/2)$, with
\begin{itemize} \item $v^2\in 4\Z$,
\item $v\cdot N_i\in 2\mathbb{Z}$  ($v\neq 2\hat{N}$, i.e. $v\in
N$ but $v/2\notin N$), \item $L^2\equiv -v^2\ \ \mod8$.
\end{itemize}}

\smallskip

In the  Proposition \ref{proposition: uniqueness of overlattice}
we prove the unicity (up to isometry) of $\mathcal{L}_{2d}'$, then
we describe the coarse moduli space of K3 surfaces with an even
set of eight disjoint rational curves (Corollary \ref{moduli}).
Moreover  we describe some known results on the N\'eron Severi
lattices of algebraic K3 surfaces with a Nikulin involution. Using
the results of \cite{bertio} we describe the relation  between the
N\'eron Severi group of an algebraic K3 surface $Y$ admitting a
Nikulin involution $\iota$ and the N\'eron Severi group of a K3
surface admitting an even set, which is the desingularization of
$Y/\iota$ (Corollary \ref{corollary: NS of X and of Y}). In the
Section \ref{section: ampleness} we analyze the ampleness of some
divisors (or more in general the nefness). These classes are used
in the Section \ref{section:projective models} to describe
projective models of algebraic K3 surfaces with an even set of
eight disjoint rational curves. In particular we describe the
following projective models: \begin{itemize} \item \textbf{double
covers of $\mathbb{P}^2$}: these branch along a sextic with eight
nodes (Paragraph \ref{cl2}) or along a smooth sextic (Paragraph
\ref{cl6}, a)) (these two situations are studied also by Barth in
\cite{barth2}, first and second cases), or along a sextic with
four nodes (Paragraph \ref{cl10}); \item\textbf{quartic surfaces
in $\mathbb{P}^3$}: these have an even set of nodes (Paragraph
\ref{cl4}) or an even set of lines (Paragraph \ref{cl8}, a))
(these two situations are studied also by Barth in \cite{barth2},
third and forth cases), or it has a mixed even set of nodes and
conics (Paragraph \ref{cl12}, b));\item \textbf{double covers of a
cone :} these branch along a conic and a sextic on the cone, which
intersect in six points (Paragraph \ref{cl4'}); \item
\textbf{complete intersections of a hyperquadric and a cubic
hypersurface in $\mathbb{P}^4$}: these have an even set of nodes
(Paragraph \ref{cl6}, b)) or
 an even set of lines (Paragraph \ref{cl10}, a));
\item \textbf{complete intersections of three hyperquadrics in
$\mathbb{P}^5$}: these have an even set of nodes (Paragraph
\ref{cl8}, b) and Paragraph \ref{cl8'}, b)) or an even set of
lines (Paragraph \ref{cl12}, a) and Paragraph \ref{cl12'}, a));
\item \textbf{double covers of a smooth quadric}: these branch
along a curve of bidegree $(4,4)$ (Paragraph \ref{cl8'}) (this
case is studied also by Barth in \cite{barth2}, sixth case); \item
We study also the following complete intersections (c.i.) of
hypersurfaces of bidegree $(a,b)$ in $\PP^n\times \PP^m$:
$$
\begin{array}{c|c|c}
\mbox{space}&\mbox{c.i.}& \mbox{paragraph}\\
\hline
\PP^1\times \PP^2&(2,3)& \ref{cl12'}~b)\\
\PP^4\times \PP^2&(2,0),(1,1),(1,1),(1,1)&\ref{cl6}~c)\\
\PP^2\times \PP^2&(1,2), (2,1)&\ref{cl16'}\\
\PP^3\times \PP^3&(1,1),(1,1),(1,1),(1,1)&\ref{cl24'}
\end{array}
$$
\end{itemize}
In Section \ref{section:projective models}  we describe moreover geometric properties of these K3 surfaces with an even set. In Section \ref{geometric} we use these
properties to give sufficient conditions for a K3 surface to have an even set.\\
{\it We would like to thank Bert van Geemen for his encouragements
and for many useful and very interesting discussions. This work
has been done during the second author stay at the University of
Milan, she would like to express her thanks to Elisabetta Colombo
 and Bert van Geemen for their warm hospitality.}


\section{K3 surfaces with an even set of nodes and of rational curves}\label{section:evenset}
\begin{defi}
Let $X$ be a surface. A set of $m$ disjoint
$(-2)$-rational smooth curves, $N_1,$ $\ldots$,$N_m$, on $X$, is an {\rm even set
of rational curves} if there is a divisor $\delta\in Pic(X)$ such
that
\begin{eqnarray*}
N_1+\ldots+N_m\sim 2 \delta,
\end{eqnarray*}
where $''\sim''$ denotes linear equivalence.
\end{defi}
\begin{defi}
Let $\bar{X}$ be a surface and let
$\mathcal{N}=\{p_1,\ldots,p_m\}$ be a set of nodes on $\bar{X}$.
Let $\tilde{\beta}:X\longrightarrow \bar{X}$ be the minimal
resolution of the nodes of $X$ and let
$N_i=\tilde{\beta}^{-1}(p_i)$, $i=1,\ldots, m$. These are $(-2)$-rational
curves on $X$. The set $\mathcal{N}$ is an {\rm even set of nodes}
if $N_1,\ldots,N_m$ are an even set of rational curves.
\end{defi}
In the case of $K3$ surfaces linear equivalence is the same as
algebraic equivalence (which we denote by $\equiv$) and
$Pic(X)=NS(X)$.\\
The existence of an even set $N_1,\ldots,N_m$ on a surface $X$ is
equivalent to the existence of a double cover
$\pi:\widetilde{Y}\rightarrow X$ from a surface $\widetilde{Y}$ to $X$
branched on $N_1+\ldots+N_m$ \cite[Lemma 17.1]{bpv}.\\

Let $Y$ be a surface and $\iota$ be an involution on $Y$ with
exactly $m$ distinct fixed points $q_1,\ldots,q_m$ and let $\widetilde{Y}$ be the
blow up of $Y$ at the points $q_1,\ldots,q_m$. The involution
$\iota$ induces an involution $\widetilde{\iota}$ on
$\widetilde{Y}$. Let $\bar{X}$ be the quotient surface $Y/\iota$
and $\pi':Y\rightarrow \bar{X}$ be the projection.
The surface $\bar{X}$ has $m$ nodes in $\pi'(q_i)$, $i=1,\ldots,
m$. Let $\widetilde{\beta}:X\rightarrow\bar{X}$ be the minimal resolution
of $\bar{X}$. Then the following diagram commutes
\begin{eqnarray}\label{diagrammone}
\begin{array}{ccc} \widetilde{Y}&\stackrel{\beta}{\lra}&Y\\
\pi \downarrow&&\downarrow \pi'\\
X&\stackrel{\tilde{\beta}}\lra& \bar{X}.
\end{array}
\end{eqnarray}
The double cover $\pi:\widetilde{Y}\rightarrow X$ is branched
on $N_1+\ldots+N_m$ where $N_i$ are the $(-2)$-curves such that
$\widetilde{\beta}(N_i)=\pi'(q_i)$, $i=1,\ldots,m$ and these form an even set.\\
Conversely if $\pi:\widetilde{Y}\rightarrow X$ is a double cover
of $X$ branched on the divisor $N_1+\ldots+N_m$ where $N_i$ are
$(-2)$-rational curves, then there is a diagram as
\eqref{diagrammone}.\\

We recall some facts about even sets on K3 surfaces:
\begin{itemize}
\item If $N_1,\ldots,N_m$ is an even set of disjoint curves on a
K3 surface, by a result of Nikulin \cite[Lemma 3]{nik} we have
$m=0,8$ or $16$.\item If $m=16$ the surface $Y$ in the diagram
\eqref{diagrammone} is a torus of dimension two (\cite[Theorem
1]{nik}), the involution $\iota$ is defined on $Y$ as $y\mapsto
-y$, $y\in Y$ and has sixteen fixed points. So $X$ is the Kummer
surface associated to the surface $Y$ (a Kummer surface is by
definition the K3 surface obtained as the desingularization of the
quotient of a torus $Y$ by the involution $y\mapsto -y$, $y\in
Y$). If $Y$ is an algebraic torus (so an Abelian surface), then
$X$ is an algebraic K3 surface and its Picard number is $\rho\geq
17$. \item If $m=8$ then the surface $Y$ is a K3 surface and the
cover involution has eight isolated fixed points (it is a {\it Nikulin
involution}, cf. Definition \ref{nik} below). If $Y$ is an
algebraic K3 surface, then $X$ is algebraic and its Picard number
is $\rho\geq 9$.
\end{itemize}

\begin{defi}\label{nik}
Let $Y$ be a K3 surface. Let $\iota$ be an involution of $Y$. The
involution $\iota$ is called {\rm Nikulin involution} if
$\iota_{|H^{2,0}(X,\C)}=id_{|H^{2,0}(X,\C)}$.
\end{defi}
We recall some facts:
\begin{itemize}
\item An involution $\iota$ on a K3 surface is a Nikulin
involution if and only if it has eight isolated fixed points
\cite[Section 5]{Nikulin symplectic}. \item The Nikulin
involutions are the unique involutions on a K3 surface $Y$ such
that the desingularization of $Y/\iota=\bar{X}$ is a K3 surface.
In fact let $\widetilde{Y}$ be the blow up of $Y$ on the fixed
points of the involution $\iota$. In this way we obtain more
algebraic classes on $\widetilde{Y}$, but the transcendental
classes are the same and so
$H^{2,0}(Y)^{\iota^*}=H^{2,0}(\widetilde{Y})^{\tilde{\iota}^*}$.
Since an automorphism of a K3 surface induces a Hodge isometry on
the second cohomology group we have
$\iota^*(H^{2,0}(Y))=H^{2,0}(Y)\simeq\C$ and since
$H^{2,0}(\widetilde{Y})^{\tilde{\iota}^*}=H^{2,0}(X)\simeq\C$ it
follows that $\iota^*$ is the identity on $H^{2,0}(Y)$, so $\iota$
is a symplectic automorphism.
\end{itemize}
\section{Even sets and Nikulin involutions}\label{section:evennik}
Let $N_1,\ldots,N_8$ be an even set of eight disjoint smooth
rational curves on a K3 surface $X$, then by adjunction $N_i^2=-2$
and Morrison shows in \cite[Lemma 5.4]{morrison} that the minimal
primitive sublattice of $H^2(X,\Z)$ containing these $(-2)$-curves
is isomorphic to the Nikulin lattice:
\begin{defi}\cite[Definition 5.3]{morrison}\label{latticenik}
The {\rm Nikulin lattice} is an even lattice $N$ of rank eight
generated by $\{N_i\}_{i=1}^8$ and $\hat N=\frac{1}{2}\sum N_i$,
with bilinear form induced by
\begin{eqnarray*}
N_i\cdot N_j=-2\delta_{ij}.
\end{eqnarray*}
\end{defi}
Observe that $\hat{N}^2=-4$ and $\hat{N}\cdot N_i=-1$. This
lattice is a negative definite lattice of discriminant $2^6$ and discriminant group $(\Z/2\Z)^{\oplus 6}$.\\
From now on $X$ is an algebraic K3 surface. A K3 surface has an
even set of eight disjoint rational curves if there are eight
disjoint rational curves spanning a copy of $N$ in $NS(X)$ (then
rank $NS(X)\geq 8$). Since $X$ is algebraic the signature of the
N\'eron Severi group $NS(X)$ is $(1,\rho-1)$, where $\rho$ is the
Picard number of $X$ (i.e. the rank of $NS(X)$). So the N\'eron
Severi group of $X$ has signature $(1,\rho-1)$ and has to contain
the negative lattice $N$ of rank eight, so $NS(X)$ contains also a
class with positive self intersection. Clearly $\rho\geq 9$ and we
will see that the generic algebraic K3 surface with an even set
has $\rho=9$ and that the number of moduli is $20-9=11$ (Corollary
\ref{moduli}). Here we study the case of algebraic K3 surfaces
with Picard number nine.
\begin{prop}\label{proposition: possible NS(X)}
Let $X$ be an algebraic K3 surface with an even set of eight
disjoint rational curves and with Picard number nine, let $L$ be a
divisor generating  $N^{\perp}\subset NS(X)$, $L^2>0$. Let $d$ be
a positive integer such that $L^2=2d$ and let
\begin{eqnarray*}
\mathcal{L}_{2d}=\Z L\oplus N.
\end{eqnarray*}
Then\\
(1) if $L^2\equiv 2$$\mod 4$ then $NS(X)=\mathcal{L}_{2d}$,\\
(2) if $L^2\equiv 0$$\mod 4$ then either $NS(X)=\mathcal{L}_{2d}$
or $NS(X)=\mathcal{L}_{2d}'$, where $\cl_{2d}'$ is generated by
$\cl_{2d}$ and by a class $(L/2,v/2)$, with
\begin{itemize} \item $v^2\in 4\Z$,
\item $v\cdot N_i\in 2\mathbb{Z}$  ($v\neq 2\hat{N}$, i.e. $v\in
N$ but $v/2\notin N$), \item $L^2\equiv -v^2\ \ \mod8$.
\end{itemize}
\end{prop}
\bprf
The discriminant group of $\cl_{2d}=\Z L\oplus N$ is
$(\Z/2d\Z)\oplus(\Z/2\Z)^{\oplus 6}$, hence an element in the
N\'eron Severi group of $X$ but not in $\cl_{2d}$ is of the form
$(\alpha L/2d,v/2)$ with $\alpha\in \Z$, $v\in N$. Since $2\cdot (\alpha
L/2d,v/2)-v\in NS(X)$ we can assume that $\alpha=d$ and so the
element is $(L/2,v/2)$. We can write $v=\sum \alpha_i
N_i+\beta\hat{N}$, $\beta\in\{0,1\}$, we have
\begin{eqnarray*}
(\frac{L}{2},\frac{v}{2})\cdot N_i\in \Z.
\end{eqnarray*}
Hence by doing the computations it follows
\begin{eqnarray*}
\frac{1}{2}(-2\alpha_i-\beta)\in \Z,
\end{eqnarray*}
hence $\beta\in 2\Z$, and so we may assume $\beta=0$. We have also
\begin{eqnarray*}
(\frac{L}{2},\frac{v}{2})\cdot \hat{N}\in \Z
\end{eqnarray*}
so
\begin{eqnarray*}
-\frac{1}{2}(\sum \alpha_i)\in\Z
\end{eqnarray*}
hence $\alpha_1+\ldots+\alpha_8\in2\Z$ and so
$\alpha_1^2+\ldots+\alpha_8^2\in 2\Z$ too.
 We have
\begin{eqnarray*}
\begin{array}{lll}
v^2&=&-2\sum\alpha_i^2-4\beta^2-2\beta\sum\alpha_i\\
&=&-2(\sum \alpha_i^2).
\end{array}
\end{eqnarray*}
It follows that $v^2\in 4\Z$ and $v\cdot N_i\in 2\mathbb{Z}$.\\
Since the N\'eron Severi lattice of a K3 surface is even we have
\begin{eqnarray*}
(\frac{L}{2},\frac{v}{2})^2=\frac{L^2+v^2}{4}\in 2\Z
\end{eqnarray*}
which gives $L^2\in 4\Z$, so $d$ must be even and $L^2+v^2\equiv
0\ \mod 8$.\\
Assume now that there is another class $(L/2,v'/2)$$\in NS(X)$,
then the class $(L/2,v/2)-(L/2,v'/2)=(v-v')/2\in NS(X)$ too. Since
$N$ is primitive $(v-v')/2\in N$. So there is a $\delta\in N$ s.t.
$v-v'=2\delta$. So $(L/2,v'/2)\in NS(X)$ if and only if
$(L/2,v'/2)=(L/2,v/2)+\delta$ for certain $\delta\in N$.
This
concludes the proof of the proposition.\eprf
\begin{prop}\label{proposition: uniqueness of overlattice}
Under the assumptions of the Proposition \ref{proposition:
possible NS(X)}, $\cl'_{2d}$ is the unique even lattice (up to isometry) such that
$[\cl'_{2d}:\cl_{2d}]=2$ and $N$ is a primitive sublattice of
$\cl'_{2d}$.
\end{prop}
\bprf  We describe briefly the group $O(N)$ of isometries of $N$.
These  must preserve the intersection form, so the image of each
$(-2)$-vector under an isometry is a $(-2)$-vector. The only
$(-2)$-vectors in the Nikulin lattice $N$ up to the sign are the
eight vectors $N_i$ and so if $\sigma\in O(N)$ then
$\sigma(N_i)=\pm N_j$, $i,j=1,\ldots,8$. In particular the group
of
permutation of eight elements $\Sigma_8$ is contained in $O(N)$. This group fixes the class $\hat{N}$.\\
Each class $v$ in $N$ is $v=\sum_{i=1}^8\alpha_i N_i+a\hat{N}$, $\alpha_i,\
a\in \Z$. We consider two different elements $v$ and $v'$ such that $\cl_{2d}$
together with the class $(L/2,v/2)$ or with the class $(L/2,v'/2)$ generate an
overlattice of $\mathcal{L}_{2d}$. We want to prove that there exists an
isometry $\sigma$ of $N$ such that $\sigma(v)=v'$. From the conditions given on
$v$, or $v'$, (in particular from the fact that $v\cdot N_i\in 2\Z$),
$v=\sum_{i=1}^8\alpha_i N_i$, $\alpha_i\in \mathbb{Z}$ and
$v'=\sum_{i=1}^8\beta_i N_i$, $\beta_i\in \Z$, we may
 assume that $\alpha_i,\beta_i\in\{0,1\}$.
The only possibilities for $v^2$ (or $v'^2$) are $-4$,$-8$,$-12$,
(by the condition on $v^2$ given in the previous proof) and this
depends only on the number of
$\alpha_i's$ (resp. $\beta_i's$) equal to one.\\
We distinguish two different cases: $v^2=v'^2$
and $v^2\neq v'^2$ but $v^2\equiv v'^2\ \mod 8$ (since $-v'^2\equiv L^2\equiv-v^2\ \mod 8$).\\
\textit{The case $v^2=v'^2$.} This condition implies that there
are the same number of $\alpha_i$ and $\beta_i$ equal to one. 
Hence there is a permutation $\sigma\in\Sigma_8\subset O(N)$ of the $N_i$, s.t. $\sigma(v)=v'$.\\
Observe that if $L^2\equiv 0\ \mod 8$, then it is clear
from the description above that $v^2=v'^2=-8$ and so we are in this case.\\
\textit{The case $v^2\neq v'^2$, $v^2\equiv v'^2 \mod 8$ and
$L^2\equiv 4$ $\mod 8$.} If $L^2\equiv 4\ \mod 8$ then $v^2$ and
$v'^2$ are $-4$ or $-12$. So we can assume $v^2=-4$ and $v'^2=-12$
and $v=N_1+N_2$ $v'=N_3+N_4+N_5+N_6+N_7+N_8$ (up to isometry
of the lattice). Observe that $v'/2=\hat{N}-v/2$ hence the lattice
generated by $\cl_{2d}$ and by $(L/2,v/2)$ or by $\cl_{2d}$ and by
$(L/2,v'/2)$ are the same.
\eprf
\begin{cor}\label{corollary: two divisible classes} Let $L^2\equiv0$$\mod 4$ and $NS(X)=\mathcal{L}_{2d}'$. Then
there are two possibilities:
\begin{itemize}
\item $L^2\equiv 4\ \mod 8$. In this case one can assume that
$v=-N_1-N_2$ and
$(L-N_3-\ldots-N_8)/2=(L+v)/2+\hat{N}-(N_3+\dots+N_8)$ is in
$NS(X)$ too.
\item $L^2\equiv 0\ \mod 8$.  In this case one can assume
that $v=-(N_1+N_2+N_3+N_4)$ and $(L-N_5-N_6-N_7-N_8)/2$ is in
$NS(X)$ too.
\end{itemize}
\end{cor}
\begin{prop}\label{suriettiva}
Let $\Gamma=\cl_{2d}$ or $\cl_{2d}'$ then there exists a K3 surface $X$ with an
even set of eight disjoint rational $(-2)$-smooth curves, such that
$NS(X)=\Gamma$.
\end{prop}
In the proof of this proposition we will use the relations between
the N\'eron Severi group of a K3 surface $Y$ with a Nikulin involution
and the N\'eron Severi group of a K3 surface $X$ which is the
desingularization of the quotient of $Y$ by the Nikulin
involution. Here we recall the two following Propositions of
\cite{bertio} in which the properties of the N\'eron Severi group
of a K3 surface with a Nikulin involution are described (we
use the notation of the Diagram \ref{diagrammone}).\\
{\bf Proposition \cite[Proposition 2.2]{bertio}} {\it Let $Y$ be an algebraic K3 surface
admitting a Nikulin involution and with Picard number nine. Let $M$ be a
divisor generating $E_8(-2)^{\perp}\subset NS(Y)$, $M^2=2d'>0$ and let
\begin{eqnarray*}
\mathcal{M}_{2d'}=\Z M\oplus E_8(-2).
\end{eqnarray*}
Then $M$ is ample, and \\
(1) if $M^2\equiv 2$$\mod 4$ then $NS(Y)=\mathcal{M}_{2d'}$,\\
(2) if $M^2\equiv 0$$\mod 4$ then either $NS(Y)=\mathcal{M}_{2d}$ or
$NS(Y)=\mathcal{M}_{2d'}'$, where $\mathcal{M}_{2d'}'$ is generated by
$\mathcal{M}_{2d}$ and by a class $(L/2,v/2)$, with $v\in E_8(-2)$.}\\
{\bf Proposition \cite[Proposition 2.7]{bertio}} {\it (1) Assume
that $NS(Y)=\Z M\oplus E_8(-2)=\mathcal{M}_{2d'}$. Let
$E_1,\ldots,E_8$ be the exceptional divisors on $\tilde{Y}$. Then:
(i) In case $M^2=4n+2$, there exist line bundles $L_1,L_2\in
NS(X)$ such that for a suitable numbering of these $E_i$ we
have:\\
$\beta^*M-E_1-E_2=\pi^*L_1,\qquad
\beta^*M-E_3-\ldots-E_8=\pi^*L_2.$ \\
The decomposition of $H^0(Y,M)$ into $\iota^*$-eigenspaces is:\\
$H^0(Y,M)\cong \pi^*H^0(X,L_1)\oplus \pi^*H^0(X,L_2),\qquad
(h^0(L_1)=n+2,\;h^0(L_2)=n+1)$\\ and the eigenspaces $\PP^{n+1},
\PP^{n}$ contain six, respectively two, fixed points.\\
(ii) In case $M^2=4n$, for a suitable numbering of the $E_i$ we have:\\ $
\beta^*M-E_1-E_2-E_3-E_4=\pi^*L_1,\qquad \beta^*M-E_5-E_6-E_7-E_8=\pi^*L_2$
with $L_1,L_2\in NS(X)$. The decomposition of $H^0(Y,M)$ into
$\iota^*$-eigenspaces is:\\ $ H^0(Y,M)\cong \pi^*H^0(X,L_1)\oplus
\pi^*H^0(X,L_2),\qquad
(h^0(L_1)=h^0(L_2)=n+1). $\\ and each of the eigenspaces $\PP^n$
contains four fixed points.\\
(2) Assume $NS(Y)=\mathcal{M}_{2d'}'$.
Then there is a line bundle $L\in NS(X)$ such that:\\
$ \beta^*M\cong\pi^*L$. The decomposition of $H^0(Y,M)$ into
$\iota^*$-eigenspaces is:\\ $H^0(Y,M)\cong H^0(X,L)\oplus
H^0(X,L-\hat{N}),\qquad (h^0(L)=n+2$, $h^0(L-\hat{N})=n)$\\ and
all fixed points map to the
eigenspace $\PP^{n+1}\subset \PP^{2n+1}$}.\\

\smallskip

{\it Proof of Proposition \ref{suriettiva}.} First observe that the lattices $\cl_{2d}$ and $\cl'_{2d}$ are primitively
embedded in the K3 lattice by \cite[Theorem 1.14.1]{Nikulin bilinear}, so we can identify them with sublattices of $U^3\oplus E_8(-1)^2$.\\
1) We consider first the case of $\Gamma=\cl_{2d}=\Z L\oplus N$ in
this case $L^2\equiv 2 \mod 4$ or $L^2\equiv 0 \mod 4$. We show
that there exists a K3 surface with an even set of $(-2)$-smooth
curves s.t. $NS(X)=\cl_{2d}$. Let $Y$ be a K3 surface with
$\rho(Y)=9$, with Nikulin involution and N\'eron Severi group of
index two in the lattice $\Z M\oplus E_8(-2)$, with $M^2\equiv 0
\mod 4$, such a K3 surface exists by \cite[Proposition 2.2,
2.3]{bertio}. We have a diagramm like Diagramm \ref{diagrammone},
and so a K3 surface $X$, which is the minimal resolution of the
quotient of $Y$ by the Nikulin involution. Since $\rho(Y)=9$ then
$\rho(X)=9$ too. By \cite[Proposition 2.7]{bertio} there is a line
bundle $L$, $L\in NS(X)$ with $\pi^*L=\beta^* M$. By the
properties of the map $\pi^*$, $2L^2=(\pi^*L)^2=(\beta^*
M)^2=M^2\equiv 0 \mod 4$ and so $L^2\equiv 2 \mod 4$ or $L^2\equiv
0 \mod 4$. Moreover $X$ has an even set made up by the eight
curves in the resolutions of the nodes of the quotient
$\bar{X}$.\\
If $L^2\equiv 2\mod 4$ then by the Proposition \ref{proposition:
possible NS(X)} $NS(X)=\cl_{2d}$, where $L^2=2d$, as
required.\\
If $L^2\equiv 0\mod 4$ we must exclude that $NS(X)=\cl_{2d}'$.
Assume that we have an element $L_1=(L-N_1-N_2)/2\in NS(X)$. We
use now the proof of \cite[Proposition 2.7]{bertio}. If
$NS(Y)=\mathcal{M}_{2d'}'$ the primitive embedding of $NS(Y)$ in
$U^3\oplus E_8(-1)^2$ is unique up to isometry. Assume that
$M^2=4n$ and choose an $\alpha\in E_8(-1)$ with $\alpha^2=-2$ if
$n$ is odd and $\alpha^2=-4$ if $n$ is even. Let $v\in
E_8(-2)\subset U^3\oplus E_8(-1)^2$ be $v=(0,\alpha,-\alpha)$ and
let $M$ be $M=(2u,\alpha,\alpha)\in U^3\oplus E_8(-1)^2$ where
$u=e_1+\frac{(n+1)}{2} f_1$ if $n$ is odd, and
$u=e_1+(\frac{n}{2}+1) f_1$ if $n$ is even (here $e_1,f_1$ denotes
the standard basis of the first copy of $U$). Then $M^2=4n$ and
$(M+v)/2=(u,\alpha,0)\in  U^3\oplus E_8(-1)^2$. This gives a
primitive embedding of $NS(Y)$ in $U^3\oplus E_8(-1)^2$, which
extends the standard one of $E_8(-2)\subset U^3\oplus E_8(-1)^2 $.
Now we can assume that $L=(u,0,\alpha)\in U(2)\oplus N\oplus
E_8(-1)\subset H^2(X,\Z)$, so by \cite[Proposition 1.8]{bertio} we
have $\beta^* M=\pi^* L$. Now
$(L-N_1-N_2)/2=(u,-N_1-N_2,\alpha)/2\in NS(X)$. By using \cite[
Proposition 1.8]{bertio} again we obtain
$\pi^*((L-N_1-N_2)/2)=(u,\frac{\alpha}{2},\frac{\alpha}{2},
-E_1-E_2)\in NS(\tilde{Y})$ and so $
(u,\frac{\alpha}{2},\frac{\alpha}{2})\in NS(Y)$, this means that
$M/2\in NS(Y)$ which is not the case. Hence $(L-N_1-N_2)/2\notin
NS(X)$, in a similar way one shows that
$(L-N_1-N_2-N_3-N_4)/2\notin
NS(X)$ and so we conclude that $NS(X)=\cl_{2d}$.\\
2) Assume now that $\Gamma=\cl'_{2d}$. In this case  we have
either\\
a) $L^2\equiv 4 \mod 8$ and so $(L-N_1-N_2)/2$ and
$(L-N_3-\ldots-N_8)/2$ are in $\Gamma$ or \\
b) $L^2\equiv 0\mod 8$ and so $(L-N_1-N_2-N_3-N_4)/2$ and
$(L-N_5-N_6-N_7-N_8)/2$ are in $\Gamma$. We do the proof assuming
that we are in case a), for the case b) the proof is very
similar.\\ Let $Y$ be a K3 surface with $\rho(Y)=9$, Nikulin
involution , N\'eron Severi group $NS(Y)=\Z M\oplus E_8(-2)$ and
$M^2=4n+2$, such a K3 surface exists by \cite[Proposition 2.2,
2.3]{bertio}. Moreover by \cite[Proposition 2.7]{bertio} there are
line bundles $L_1$ and $L_2$ in $NS(X)$ with $\beta^*
M-E_1-E_2=\pi^* L_1$, $\beta^*M-E_3-\ldots-E_8=\pi^* L_2$. Since
the embedding of $\Z L\oplus E_8(-2)$ in the K3 lattice is unique
we may assume that $M=e_1+(2n+1)f_1$ and
$L_1=(e_1+(2n+1)f_1+N_1+N_2)/2-N_1-N_2\in NS(X)$, by
\cite[Proposition 1.8]{bertio} we have $\beta^*M-E_1-E_2=\pi^*
L_1$. The class $U(2)\ni (e_1+(2n+1)f_1)=2 L_1+N_1+N_2$ is in
$NS(X)$, is orthogonal to the $N_i$ and has self intersection
$8n'+4$, we call it $L$. By Proposition \ref{proposition: possible
NS(X)} we have $NS(X)=\cl_{2d}'$ with $d=4n+2$, so we are done.
\eprf 
\textbf{Remark.} By using the {\it surjectivity of the
period map} one can show the existence of a K3 surface $X$ with
$NS(X)=\Gamma$, it is however difficult to show that there is an
embedding of the classes $N_i$ as irreducible $(-2)$-smooth curves
in $NS(X)$. This is assured by the previous proposition.\\
From the Proposition \ref{suriettiva} follows a relation between
the N\'eron Severi group of the K3 surface $Y$ admitting a Nikulin
involution and the N\'eron Severi group of a K3 surface $X$ which
is the desingularization of the quotient.
\begin{cor}\label{corollary: NS of X and of Y} Let $Y$ be an
algebraic K3 surface with $\rho(Y)=9$ admitting a Nikulin
involution, and let $X$ be the desingularization of its quotient.\\
(1) $NS(Y)=\mathcal{M}_{2d}$ if and only if $NS(X)=\mathcal{L}_{4d}'$;\\
(2) $NS(Y)=\mathcal{M}_{4d}'$ if and only if
$NS(X)=\mathcal{L}_{2d}$.\end{cor} \bprf The proof follows from
\cite[Proposition 2.7]{bertio} and Proposition
\ref{suriettiva}. We sketch it briefly.\\
The proof of the direction $\Leftarrow$ of the statement follows immediately from the proof of Proposition \ref{suriettiva}. For the other direction
we distinguish three cases (we use the notation of \textit{loc. cit.}):\\
(a) Case (1), $(i)$. Clearly $(\beta^*M-E_1-E_2)^2=(\pi^*L_1)^2$
and in the proofs of Proposition \ref{suriettiva}, case (2), and
of \cite[Proposition 2.7]{bertio} it is proved that
\begin{eqnarray}\label{formula: L1 and L2: first case } L_1=(L-N_1-N_2)/2,\ \
L_2=(L-N_3-\ldots-N_8)/2.\end{eqnarray} Since  $\pi$ is a $2:1$
map to $X$ the previous equality becomes
$4n+2-1-1=\frac{1}{2}(L-N_1-N_2)^2=\frac{1}{2}(L^2-4)$ and so
$L^2=2(4n+2)$. By the Proposition \ref{proposition: possible NS(X)}, where we describe the possible N\'eron-Severi groups of K3 surfaces with an even set, we obtain that  $NS(X)=\cl'_{2d}$, $d\equiv 2 \mod 4$.\\
(b) Case (1), $(ii)$. As before, in the proof of \cite[Proposition
2.7]{bertio} it is proved that:
\begin{eqnarray}\label{formula: L1 and L2: second case }
L_1=(L-N_1-\ldots-N_4)/2,\ L_2=(L-N_5-\ldots-N_8)/2.
\end{eqnarray}
So we obtain
$4n-4=(\beta^*M-E_1-E_2-E_3-E_4)^2=2((L-N_1-\ldots-N_4)/2)^2=
\frac{1}{2}(L^2-8)$ and so $L^2=2(4n)$. By the Proposition
\ref{proposition: possible NS(X)} we obtain that
$NS(X)=\cl'_{2d}$, $d\equiv 0\mod
4$.\\
(c) Case (2), $M^2=2L^2$, and so by an argumentation as in the
proof of the Proposition \ref{suriettiva}, case (1), we have
$NS(X)=\mathcal{L}_{2d}$.\eprf

Some explicit correspondence between the K3 surfaces $Y$ and $X$ are shown in the Table \ref{tab1}.\\
\textbf{Remark.} Let $X$ be a K3 surface such that the lattice
$\Gamma=\cl_{2d}$ or $\cl_{2d}'$ is primitively embedded in $NS(X)$
 and $\rho(X)\geq 9$. There exists a
deformation of the K3 surface $\{X_t\}$ such that
$X_{\overline{t}}=X$ and $X_0$ is such that $NS(X_0)=\Gamma$. Let
$Y_0$ be the K3 surface such that the desingularization of its
quotient by a Nikulin involution is $X_0$. The N\'eron Severi
group of $Y_0$ is either $\mathcal{M}_{4d}'$ or $\mathcal{M}_{d}$.
The deformation on $X$ induces a deformation $\{Y_t\}$ of $Y_0$ such that
 the surface $Y_{\overline{t}}$ admits a Nikulin involution and the
desingularization of its quotient by the Nikulin involution is
$X_{\overline{t}}$. This means that $X_{\overline{t}}$ admits
an even set of eight disjoint rational curves.\\
In particular if $X$ is an algebraic K3 surface such that
$\cl_{2d}$ (resp. $\cl_{2d}'$) is primitively embedded in $NS(X)$,
then $X$ is the minimal resolution of the quotient of a K3 surface
$Y$ such that
$\mathcal{M}_{4d}'$ (resp. $\mathcal{M}_{d}$) is primitively
embedded in $NS(Y)$.

\begin{cor}\label{moduli}
The coarse moduli space of $\Gamma$-polarized $K3$ surfaces (cf.
\cite[p.5]{dolgachev} for the definition) is the quotient of
\begin{eqnarray*}
\mathcal{D}_{\Gamma}=\{\omega\in\PP(\Gamma^{\perp}\otimes_{\Z}\C):~\omega^2=0,~\omega\bar{\omega}>0\}
\end{eqnarray*}
by an arithmetic group $O(\Gamma)$ and has dimension eleven. The generic K3 surface with an even set of eight disjoint rational curves
has Picard number nine.\\
\end{cor}
\bprf By Proposition \ref{proposition: possible NS(X)}, each K3
surface with an even set is contained in this space, on the other
hand, by Proposition \ref{suriettiva} each point of this space
corresponds to a K3 surface with an even set of irreducible
$(-2)$-curves. Moreover the generic K3 surface in this space has
Picard number nine. 


By using the results on lattices of K3 surfaces with a Nikulin
involution and with an even set (\cite[Proposition 2.2]{bertio}
and Proposition \ref{proposition: possible NS(X)}) it is possible
to prove that certain K3 surfaces admitting a Nikulin involution
do not admit an even set and viceversa.

\begin{lemma}\label{nikotto}
Let $\mathcal{L}_{2d}'$, $\mathcal{L}_{2d}$, $\mathcal{M}_{2d'}'$,
$\mathcal{M}_{2d'}$ be the lattices described in the Proposition
\ref{proposition: possible NS(X)} and in \cite[Proposition 2.2]{bertio}. Then the discriminant groups of these
lattices are the following:
\begin{itemize}
\item
$(\mathcal{L}_{2d})^{\vee}/\mathcal{L}_{2d}=(\mathbb{Z}/2d\mathbb{Z})\oplus
(\mathbb{Z}/2\mathbb{Z})^{\oplus 6}$; \item
$(\mathcal{L}_{2d}')^{\vee}/\mathcal{L}_{2d}'=(\mathbb{Z}/2d\mathbb{Z})\oplus
(\mathbb{Z}/2\mathbb{Z})^{\oplus 4}$;\item
$(\mathcal{M}_{2d'})^{\vee}/\mathcal{M}_{2d'}=(\mathbb{Z}/2d'\mathbb{Z})\oplus
(\mathbb{Z}/2\mathbb{Z})^{\oplus 8}$;\item
$(\mathcal{M}_{2d'}')^{\vee}/\mathcal{M}_{2d'}'=(\mathbb{Z}/2d'\mathbb{Z})\oplus
(\mathbb{Z}/2\mathbb{Z})^{\oplus 6}$.
\end{itemize}
\end{lemma}
\bprf The discriminant group of $(\mathcal{L}_{2d})^{\vee}/\mathcal{L}_{2d}$ is
generated by the classes
$L/2d, (N_i+N_{i+1})/2$ $i=1,\ldots, 6$.\\
Now we consider the lattice $\mathcal{L}_{2d}'$. We suppose that it is
generated by $\mathcal{L}_{2d}$ and by the class $(L-N_1-N_2)/2$ (here we are
supposing that $d\equiv 2$ mod 4; if it is not so, then $d\equiv 0$ mod 4 and
the class that we need is $(L-N_1-N_2-N_3-N_4)/2$, the computation in the two
cases are essentially the same). The group
$(\mathcal{L}_{2d}')^{\vee}/\mathcal{L}_{2d}$ is generated by the classes of
$(\mathcal{L}_{2d})^{\vee}/\mathcal{L}_{2d}$ whose product with the class
$(L-N_1-N_2)/2$ is an integer. A basis for this space is $\{(L+d(N_2+N_3))/2d$,
$(N_1+N_2)/2$, $(N_3+N_4)/2$, $(N_4+N_5)/2$, $(N_5+N_6)/2$, $(N_6+N_7)/2\}$.
Now we consider $(\mathcal{L}_{2d}')^{\vee}/\mathcal{L}_{2d}'$. Since we make a
quotient with a larger lattice the discriminant group is smaller than
$(\mathcal{L}_{2d}')^{\vee}/\mathcal{L}_{2d}$, in fact $(N_1+N_2)/2\equiv
-(L-N_1-N_2/2)+d(L+d(N_2+N_3)/2d)$ $\mod \cl'_{2d}$ since $d$ is even. So a
basis for $(\mathcal{L}_{2d}')^{\vee}/\mathcal{L}_{2d}'$ is
$\{(L+d(N_2+N_3))/2d$, $(N_3+N_4)/2$, $(N_4+N_5)/2$, $(N_5+N_6)/2$,
$(N_6+N_7)/2\}$.\\
The discriminant of $(\mathcal{M}_{2d'})^{\vee}/\mathcal{M}_{2d'}$ is generated
by $M/2d'$, $E_i/2$, $i=1,\ldots, 8$ where $E_i$, $i=1,\ldots,8$
 is
the basis of the lattice $E_8(-2)$.\\
The lattice $\mathcal{M}_{2d'}'$ is generated by
$\mathcal{M}_{2d'}$ and the class $(M-E_1)/2$ or $(M-E_1-E_3)/2$.
Its discriminant group can be computed in a similar way as
before.\eprf

\begin{cor}\label{relation between NS(X) and NS(Y)} Let $X$ be a
K3 surface such that either $NS(X)=\cl_{2d}$ or $NS(X)=\cl_{2e}'$.
Then $X$ does not admit a Nikulin involution (by the Proposition
\ref{suriettiva} $X$ has an even set).\\
Let $Y$ be a K3 surface such that $NS(Y)=\mathcal{M}_{2d'}$. Then
$Y$ does not have an even set of eight disjoint rational curves (by
\cite[Proposition 2.3]{bertio} $Y$ admits a Nikulin involution).

\end{cor}
\bprf If a surface admits an even set of disjoint rational curves
and a Nikulin involution then its N\'eron Severi group is
isometric to a lattice of type $\mathcal{L}_{2d}'$ or
$\mathcal{L}_{2d}$ and to a lattice of type $\mathcal{M}_{2d'}'$
or $\mathcal{M}_{2d'}$. If two lattices are isometric, then their
discriminant group are equal. The only lattices which could have
the same discriminant group are $\mathcal{L}_{2d}$ and
$\mathcal{M}_{2d}'$ (for the same $d$), and the latter lattice
requires $d\equiv 0$ mod 4. This proves the corollary.\eprf

\section{Ampleness and nefness of some divisors on $X$}\label{section: ampleness}
Our next aim (cf. Section \ref{section:projective models}) is to
describe projective models of K3 surfaces with an even set of
eight disjoint rational curves. Here we give some results on
ampleness and on nefness of divisors on such  K3 surfaces. We
prove moreover that the associated linear systems have no base
points. These properties guaranty that the maps induced by the
linear systems are regular (in fact birational) maps.

\begin{defi}
A divisor $L$ on a surface $S$ is:
$$
\begin{array}{llll}
\bullet&\mbox{\textbf{nef }}&\mbox{if } L^2\geq 0&\mbox{and
}L\cdot C\geq 0 \mbox{ for each
irreducible curve }C\mbox{ on } S,\\
\bullet&\mbox{\textbf{pseudo ample} }&\mbox{if } L^2> 0&\mbox{and
}L\cdot C\geq 0 \mbox{ for each
irreducible curve }C\mbox{ on } S,\\
&\mbox{(or \textbf{big and nef)}}&&\\
\bullet&\mbox{\textbf{ample} }&\mbox{if } L^2> 0&\mbox{and }L\cdot
C>0 \mbox{ for each
irreducible curve }C\mbox{ on } S.\\
\end{array}
$$
\end{defi}
If $X$ is a K3 surface with a line bundle $L$ such that  $L^2\geq
0$, the condition $L\cdot C\geq 0$ for each irreducible curve $C$
on $X$ is equivalent to
the condition $L\cdot \delta\geq 0$ for each irreducible $(-2)$-curve $\delta$ on $X$ (cf. \cite[Proposition 3.7]{bpv}).\\

Let $H$ be an effective divisor on a K3 surface. The intersection of $H$ with
each curve $C$ is positive except when $C$ is a component of $H$ and $C$ is a
$(-2)$-curve. If the linear system $|H|$ does not have fixed components and if
$H^2>0$, then the generic element in $|H|$ is smooth and irreducible and $H$ is
a pseudo ample divisor (cf. \cite[Proposition 2.6]{saintdonat}). The fixed
components of a linear system on a K3 surface are always $(-2)$-curves
\cite[Paragraph 2.7.1]{saintdonat}. Recall that by \cite[Theorem p.79]{kollar}
if $H$ is pseudo ample (or ample) then either $|H|$ has no fixed components or
$H=aE+\Gamma$ where $|E|$ is a free pencil and $\Gamma$ is an irreducible
$(-2)$-curve such that $E\Gamma=1$. Finally in \cite[Corollary 3.2]{saintdonat}
Saint-Donat proves that a linear system on a K3 surface has
no base points outside its fixed components.\\

Let now $H$ be a pseudo ample divisor on $X$. If $|H|$ has a fixed
component $B$, then $H=B+M$, where $M$ is the moving part of the
linear system $|H|$. The linear system $|H|$ defines a map
$\phi_{H}$ and if $H=M+B$ then $\phi_{H}=\phi_{M}$. Now we assume that
$|H|$ has no fixed components (and hence no base points). The
system $|H|$ defines the map:
\begin{eqnarray*}
\phi_{H}:X\lra\PP^{p_a(H)}
\end{eqnarray*}
where $p_a(H)=H^2/2+1$ and there are two cases (cf. \cite[Paragraph 4.1]{saintdonat}):\\
(i) either $\phi_H$ is of degree two and its image has degree $p_a(H)-1$ ($\phi_H$ is {\it hyperelliptic}),\\
(ii) or $\phi_H$ is birational and its image has degree $2p_a(H)-2$.\\
In particular in the second case if $H$ is ample (i.e.\ does not
contract $(-2)$-curves) $\phi_{H}$ is an embedding, so $H$ is very
ample.

\begin{prop}\label{prop:pseudoample} Let $X$ be as in the Proposition \ref{proposition: possible
NS(X)}. Then we may assume that $L$ is pseudo ample and it has no
fixed components.\end{prop} \bprf By the Diagram \ref{diagrammone}
since $X$ and $Y$ are algebraic $\bar{X}$ is embedded in some
projective space and has eight nodes. The generic hyperplane
section of $\bar{X}$ is a smooth and irreducible curve (it does
not pass through the nodes). Its pull back on $X$ is then
orthogonal to $N_1,\ldots,N_8$, we call it $\mathcal{H}$, observe
that $\mathcal{H}=\alpha L$ for some integer $\alpha$. Since
$\mathcal{H}$ is pseudo ample then $L$ is pseudo ample too, in
particular observe that $L\Gamma>0$ for each $(-2)$-curve which is
not one of the $N_i$'s. If $L$ has fixed components then by
\cite[Theorem p.79]{kollar} it is $L=aE+\Gamma$ where $|E|$ is a
free pencil and $\Gamma$ an irreducible $(-2)$-curve such that
$E\Gamma=1$. If $\Gamma\not= N_i$ for each $i=1,\ldots, 8$, then
$0<L\Gamma=a-2$, which gives $a>2$. Now $0\leq LN_i= aEN_i+\Gamma
N_i$, since $\Gamma N_i\geq 0$ and $a>2$ we obtain $\Gamma N_i=0$
for each $i$, so $\Gamma$ is in $(N)^{\perp}$ which is not
possible. If $\Gamma=N_i$ for some $i$, then $0=LN_i=a-2$ so $a=2$
then $L=2E+N_i$ and so $(L-N_i)/2$ is in the N\'eron Severi group
too which is not the case. So by \cite[Proposition
2.6]{saintdonat} we can assume that $L$ is smooth and
irreducible.\eprf

\begin{prop}\label{theorem:compobase}
Let $X$ be as in the Proposition \ref{proposition: possible
NS(X)}. If $d\geq 3$, i.e. $L^2\geq 6$, then the class $L-\hat{N}$ in
the N\'eron Severi group is an ample class.
\end{prop}
\bprf The self intersection of $L-\hat{N}$ is
$(L-\hat{N})^2=2d-4$, which is positive for each $d\geq 3$. So to
prove that $L-\hat{N}$ is ample we have to prove that for each
irreducible $(-2)$-curve $C$ the intersection number
$C\cdot(L-\hat{N})$ is positive.\\
In the proof we use the inequality:
\begin{eqnarray}\label{inequality1} (\sum_{i=1}^n x_i)^2\leq n\sum_{i=1}^n x_i^2\end{eqnarray}
which is true for every $(x_1,\ldots,x_n)\in\mathbb{R}^n$. 
Suppose that there exists an effective irreducible curve
$C$ such that $C\cdot (L-\hat{N})\leq 0$, then we prove that
$C\cdot C<-2$.\\
We observe that each element in the N\'eron Severi group is a
linear combination of $L$ and $N_i$ with coefficients in
$\frac{1}{2}\Z$. We consider the curve $C=aL+\sum_{i=1}^8 b_i
N_i$ where $a,b_i\in\frac{1}{2}\Z$. If $a=0$ the only possible
$(-2)$-curves are the $N_i$'s and $N_i\cdot (L-\hat{N})=1$. So we
can assume that $a\neq 0$.
Since $C$ is an irreducible curve, it has a non-negative
intersection with all effective divisors. Hence $C\cdot L=2da\geq
0$, so $a>0$,
and $C\cdot N_i=-2b_i\geq 0$, so $b_i\leq 0$.\\
Now we assume that $C\cdot (L-\hat{N})\leq 0$, then
$$(aL+\sum_{i=1}^8 b_i N_i)\cdot (L-\hat{N})=2da+\sum_{i=1}^8 b_i\leq 0.$$
Since $b_i\leq 0$,  $2da-\sum_{i=1}^n |b_i|\leq 0$ and so $2da\leq
\sum_{i=1}^n |b_i|$, where each member is non negative. So it is
possible to pass to the square of the relation, obtaining
$4d^2a^2\leq (\sum_{i=1}^n |b_i|)^2.$ Using the relation
\eqref{inequality1} one has
\begin{eqnarray}\label{inequality2}
4d^2a^2\leq (\sum_{i=1}^n
|b_i|)^2\leq 8\sum_{i=1}^8 b_i^2.
\end{eqnarray} Now we compute
the square of $C$ and we use the inequality \eqref{inequality2}
to estimate it:
$$C\cdot C=2da^2-2\sum_{i=1}^8(b_i^2)\leq 2da^2-d^2a^2.$$
If $d\geq 5$ then $\sqrt{\frac{2}{d^2-2d}}<\frac{1}{2}$, and so
for $d\geq 5$ we have $C\cdot C<-2$ (because $a\geq\frac{1}{2}$).
This proves
the theorem in the case $d\geq 5$.\\
More in general for each $d\geq 3$, $\sqrt{\frac{2}{d^2-2d}}<1$,
so for the cases $d=3$ and $d=4$ one has to study only the case
$a=\frac{1}{2}.$\\
For $d=3$, then $L^2=6$ and so in $NS(X)$ all
the elements are of the form $aL+\sum_{i=1}^8 b_i N_i$ with
$a\in\Z$ (and not in $\frac{1}{2}\Z$). Then the theorem is proved
exactly in the same way as before.\\
Let $d=4$. The only possible irreducible
$(-2)$-curves with a negative intersection with $(L-\hat{N})$ are
of the form $\frac{1}{2} (L+N_1+N_2+N_3+N_4) +\sum_{i=1}^8 \beta_i
N_i$ with $\beta_i\in\Z$, $\beta_1,\ldots,\beta_4\leq -1$ and
$\beta_5,\ldots,\beta_8\leq 0$. The only $(-2)$-curves of this
type are $\frac{L+N_1+N_2+N_3+N_4}{2}-N_1-N_2-N_3-N_4-N_j$,
$j=5,6,7,8$  and these curves have a positive intersection with
$L-\hat{N}$. Then the proposition is proved also for $d=4$. \eprf

\begin{prop}\label{prop: mL-N e' nef}
In the situation of Proposition \ref{theorem:compobase},
$m(L-\hat{N})$ and $mL-\hat{N}$ for $m\in\Z_{>0}$, are ample. If
$d=2$, i.e. $(L-\hat{N})^2=0$, then $m(L-\hat{N})$ is nef and
$mL-\hat{N}$ is ample for $m\geq 2$.
\end{prop}
\bprf It is a similar computation as in the proof of Proposition
\ref{theorem:compobase}.\eprf

\begin{prop}\label{prop:compo1}
The divisors $L-\hat{N}$, $mL-\hat{N}$ and $m(L-\hat{N})$,
$m\in\Z_{>0}$, do not have fixed components for $d\geq 2$.
\end{prop}
\bprf We proof the proposition for the divisor $L-\hat{N}$. The
proof in the
other cases is essentially the same.\\
For $d=2$ we have $(L-\hat{N})^2=0$ and is nef by the Proposition
\ref{prop: mL-N e' nef} so by \cite[Theorem p. 79, (b)]{kollar}
$L-\hat{N}=aE$ where $|E|$ is a free pencil,
and so the assertion is proved in this case.\\
Assume $d\geq 3$, then for \cite[Theorem p. 79, (d)]{kollar} we
have either $L-\hat{N}$ has no fixed components or
$L-\hat{N}=aE+\Gamma$, where $|E|$ is a free pencil and $\Gamma$
is an irreducible $(-2)$-curve such that $E\Gamma=1$. We assume we
are in the second case, then since $L-\hat{N}$ is ample we have
$$
0<\Gamma(L-\hat{N})=a-2
$$
and so $a>2$. We distinguish two cases:\\
1. $\Gamma=\alpha L+\sum \beta_j N_j$, $\Gamma\not= N_i$ for each
$i$, so $\alpha\not=0$. For each $i$ we have:
$$
1=N_i(L-\hat{N})=aEN_i+\Gamma N_i
$$
Since $EN_i\geq 0$ and $a>2$ then $EN_i=0$ and so
$$
1=\Gamma N_i=(\alpha L+\sum_j \beta_j N_j)N_i=-2\beta_i.
$$
We obtain $\beta_j=-1/2$ for all $j$, so
$$
\Gamma=\alpha L-\frac{N_1+\ldots+N_8}{2}=\alpha L-\hat{N}.
$$
By considering the self-intersection of $\Gamma$ we obtain
$$
-2=\alpha^22d-4\geq 6\alpha^2-4
$$
which is positive since $\alpha$ is a non zero integer. So this case is not possible.\\
2. $\Gamma=N_i$ for some $i=1,\ldots, 8$. We have
$$
1=N_i(L-\hat{N})=N_i(aE+N_i)=aEN_i-2=a-2
$$
so $a=3$, $L-\hat{N}=3E+N_i$. For $j\not=i$ we have
$$
1=(L-\hat{N})N_j=3EN_j
$$
but this is impossible. Hence $L-\hat{N}$ has no base components.
\eprf

\begin{lemma}\label{lemma:modelli}
The map $\phi_{L-\hat{N}}$ is
\begin{itemize}
\item an
embedding  if  $L^2\geq 10$,
\item a 2:1 map to $\mathbb{P}^1\times\mathbb{P}^1$ if  $L^2=8$,
\item a 2:1 map to $\mathbb{P}^2$ if  $L^2=6.$
\end{itemize}
\end{lemma}
\bprf By the Proposition \ref{theorem:compobase} $L-\hat{N}$ is
ample and by the Proposition \ref{prop:compo1} $|L-\hat{N}|$ has
no fixed components; for a K3 surface this implies that
$|L-\hat{N}|$ has no base points too (cf. \cite[Corollary
3.2]{saintdonat}), and it defines a map $\phi_{L-\hat{N}}$. The
assertion for $L^2=6$ is clear since $(L-\hat{N})^2=2$ and hence
the map $\phi_{L-\hat{N}}$ defines a double cover of $\PP^2$. We
show that in the case $L^2=2d\geq 10$, i.e. $d\geq 5$, the map is
not hyperelliptic. By \cite[Theorem 5.2]{saintdonat} $L-\hat{N}$
is hyperelliptic iff $(i)$ there is an elliptic irreducible curve
$E$ with $E\cdot(L-\hat{N})=2$ or $(ii)$ there is an irreducible
curve $B$, with $p_a(B)=2$ and $L-\hat{N}=\mathcal{O}(2B)$. The
case $(ii)$ would implies $L-\hat{N}\equiv 2B$ and so
$\frac{1}{2}(L-\hat{N})\in NS(X)$ which is not possible by the
description of $NS(X)$ of Proposition \ref{proposition: possible
NS(X)}. We have to exclude $(i)$. We argue in a similar way as in
Proposition \ref{theorem:compobase}. Assume that there is $E=a
L+\sum b_iN_i$ an irreducible curve with $E\cdot (L-\hat{N})=2$.
Then we show $E^2\not=0$. Since $E$ is the class of an irreducible
curve, $a\in\frac{1}{2}\Z_{>0}$ and $b_i\in\frac{1}{2} \Z_{\leq
0}$. We have $2=E\cdot (L-\hat{N})=2da+\sum_{i=1}^{8} b_i$ and so
$2da-2=-\sum_{i=1}^{8} b_i$ which gives together with the
inequality \eqref{inequality1}:
\begin{eqnarray*}
4(da-1)^2=(\sum_{i=1}^{8}|b_i|)^2\leq 8 \sum_{i=1}^{8}|b_i|^2
\end{eqnarray*}
and so $(da-1)^2\leq 2\sum_{i=1}^{8}b_i^2$. On the other hand we
have
\begin{eqnarray*}
E^2=2da^2-2\sum_{i=1}^{8}b_i^2\leq
2da^2-(da-1)^2=2da^2-d^2a^2-1+2da.
\end{eqnarray*}
We have $E^2<0$ for $a<\frac{d-\sqrt{2d}}{d(d-2)}$ or
$a>\frac{d+\sqrt{2d}}{d(d-2)}$, since $a>1/2$ and
$\frac{d-\sqrt{2d}}{d(d-2)}<\frac{1}{2}$ for each $d\geq 5$ and
$\frac{d+\sqrt{2d}}{d(d-2)}<\frac{1}{2}$ for each $d\geq 6$, we
obtain $E^2<0$ for $d\geq 6$. We analyze the case of $d=5$. Here
$L^2=10$ and so $a\in\Z_{>0}$, for $d=5$ we have
$\frac{d+\sqrt{2d}}{d(d-2)}=\frac{5+\sqrt{10}}{15}<1$. In
conclusion for each $d\geq 5$ we obtain $E^2<0$. In the case of
$d=4$, then we have $L^2=8$ and the classes
$E_1=\frac{L-N_1-N_2-N_3-N_4}{2}$,
$E_2=\frac{L-N_5-N_5-N_7-N_8}{2}$ are in the N\'eron Severi group.
We have $E_1^2=E_2^2=0$, $E_1\cdot E_2=2$ and $L-\hat{N}=E_1+E_2$,
so $\phi_{L-\hat{N}}$ defines a 2:1 map to a quadric in
$\PP^1\times\PP^1$ (cf. \cite[Proposition 5.7]{saintdonat}). \eprf

\begin{prop}\label{prop:compobase}\label{prop:compo2}
1) Let $D$ be the divisor $D=L-(N_{1}+\ldots+N_{r})$ (up to relabel the
indices), $1\leq r \leq 8$.
\begin{itemize}
\item If $NS(X)=\cl_{2d}$, then $D$ is pseudo ample for $d>r$;
\item if $NS(X)=\cl_{2d}'$, then $D$ is nef for $d=r+4$ and pseudo
ample for $d>r+4$, \item if $D$ is pseudo ample and
$NS(X)=\cl_{2d}$ then it does not have fixed components.
\end{itemize}
2) Let $NS(X)=\cl_{2d}'$. Let
$\bar{D}=(L-(N_{1}+\ldots+N_{r}))/2$ with  $r=2,6$ if $2d\equiv 4\mod 8$ and $r=4$ if $2d\equiv 0\mod 8$.
Then
\begin{itemize}
\item the divisor $\bar{D}$ is nef and is pseudo ample whenever it has positive
self intersection, \item if $\bar{D}$ is pseudo ample then it does not have
fixed components, if $\bar{D}^2=0$ then the generic element in $|\bar{D}|$ is
an elliptic curve.
\end{itemize}
\end{prop}
\bprf The arguments are similar as those of the proof of the
Proposition \ref{theorem:compobase}.\\
1) Let $C=aL+\sum b_iN_i$ be an effective irreducible curve on
$X$, $a\in \frac{1}{2}\Z_{>0}$, $b_i\in\frac{1}{2}\Z_{\leq 0}$
such that $C\cdot (L-(N_{1}+\ldots+N_{r}))\leq 0$. We prove
that such $C$ has $C\cdot C<-2$.
If $L^2=2d$ the condition $C\cdot (L-(N_{1}+\ldots+N_{r}))\leq
0$ is equivalent to $da\leq|b_{1}|+\ldots +|b_{r}|$. By using the
inequality \eqref{inequality1} one finds
$$d^2a^2\leq (\sum_{i=1}^r b_{i})^2\leq r \sum_{i=1}^r(b_{i}^2).$$
The self intersection of the curve $C$ is $C\cdot
C=2da^2-2\sum_{i=1}^8 b_i^2$, and so
$$\begin{array}{ll}
2da^2-2\sum_{i=1}^8 b_i^2&=2da^2-2\sum_{i=1}^r
b_{i}^2-2\sum_{i=r+1,\dots,8}b_i^2\\
&\leq 2da^2-(2d^2a^2/r)-2\sum_{i=r+1,\ldots,8}b_i^2\\
&\leq 2da^2-(2d^2a^2/r)=2da^2(1-(d/r))=-2da^2((d-r)/r)\\
&\leq -2ra^2((d-r)/r)<-2.\\
\end{array}
$$
Where the last inequality holds for $d>r$ or $d>r+4$. If $d=r+4$ we have
$$
2da^2-2\sum_{i=1}^8 b_i^2\leq 2da^2(1-(d/r))=-8a^2(1+4/r)<-2.
$$
So also in this case the assertion is proved.\\
We show that for $NS(X)=\cl_{2d}$ and $D$ pseudo ample then it does not have base components. By \cite[Theorem p. 79 (d)]{kollar}, either $|D|$ does not have base components or
$D=aE+\Gamma$ where $|E|$ is a free pencil and $\Gamma$ is an irreducible $(-2)$-curve with
$E\cdot \Gamma=1$. We assume, we are in this case, moreover put $\Gamma=\alpha L+\sum_{i=1}^{8}\beta_i N_i$ with $\alpha$ a non zero integer and $\beta_i\in\Z_{\leq 0}$.\\
{\bf a)} Assume $\Gamma\not= N_i$ for each $i$. By the proof above we have $0<\Gamma D=a-2$ and so $a>2$. For $r+1\leq j\leq 8$ we have $0=DN_j=aEN_j+\Gamma N_j$ and so $EN_j=\Gamma N_j=0$. On the other hand $2=DN_1=aEN_1+\Gamma N_1$. Since $a>2$ we have $E N_1=0$ and $\Gamma N_1=2$. In fact this holds for each $N_i$, $i=1,\ldots,r$ then $E N_i=0$ for each $i=1,\ldots,8$ and so $E=\delta L$, but then $E\Gamma=2d\alpha\delta >1$, which is not possible.\\
{\bf b)} Assume $\Gamma=N_i$. {\bf $b_1$)} First we consider the case $\Gamma=N_i$, $1\leq i\leq r$, w.l.o.g $i=1$. Then $D=aE+N_1$, with $E N_1=1$. We obtain $2=DN_1=aEN_1-2$, and so $a=4$. Now $2=D N_2=4EN_2$ which is not possible. {\bf $b_2$)} Assume $i=r+1$, then $0=D N_{r+1}=a-2$ which gives $a=2$. By substituting we obtain $L-N_1-\ldots-N_{r+1}=2E$ but this is not possible by the structure of the Neron Severi group.\\
2) We show the assertion in the case $r=2$, the other cases are
similar, w.l.o.g we assume that the divisor is $(L-N_1-N_2)/2$ and
so we have $L^2=4d'$. First observe that:
\begin{eqnarray*}
\left(\frac{L-N_1-N_2}{2}\right)^2=d'-1\geq 0.
\end{eqnarray*}
Now let $C=aL+\sum b_iN_i$ be an effective  $(-2)$-curve, with $a$
and $b_i$ as before. Assume that
\begin{eqnarray}\label{eq:caso2}
0\geq C\cdot \left(\frac{L-N_1-N_2}{2}\right)=2d'a+b_1+b_2.
\end{eqnarray}
We show that this implies $C^2\not= -2$. We have
\begin{eqnarray*}
C^2=4d'a^2-2(b_1^2+b_2^2)-2\sum_{i=3}^8 b_i^2\leq
4d'a^2-(b_1+b_2)^2-2\sum_{i=3}^8 b_i^2
\end{eqnarray*}
since $2(b_1^2+b_2^2)\geq (b_1+b_2)^2$ by the inequality
\eqref{inequality1}. By using the assumption we obtain
\begin{eqnarray}\label{curve-2}
4d'a^2-(b_1+b_2)^2-2\sum_{i=3}^8 b_i^2\leq
4d'a^2-4d'^2a^2-2\sum_{i=3}^8
b_i^2=4d'a^2(1-d')-2\sum_{i=3}^{8}b_i^2.
\end{eqnarray}
We distinguish two cases:\\
{\bf 1. $d'\geq 2$}. Then we obtain with $a=a'/2$, $a'\in\Z$,
$a\not= 0$, $2(2d'a^2(1-d')-\sum_{i=3}^8 b_i^2)\leq -2a'^2$.
We have two possibilities:\\
{\bf a) $a'\geq 2$} then $-a'^2\leq -4$, and so $C^2\leq -8$.\\
{\bf b) $a'=1$}. Here we have to analyze the case  $2(2d'a^2(1-d'))=2(d'/2(1-d'))=-2$ and the inequalities in \eqref{curve-2} are equalities, which corresponds, if possible, to $C^2=-2$.\\
This gives $d'(1-d')=-2$, and so $d'=2$. Moreover we have
$-1-\sum_{i=3}^8 b_i^2=-1$ and so $b_3=\ldots=b_8=0$. Hence
$C=\frac{1}{2}L+b_1N_1+b_2N_2$ and it must be $b_1=-1/2+b_1'$,
$b_2=-1/2+b_2'$, $b_i'\in\Z$, so by using the inequality
\eqref{eq:caso2}, we obtain $d'\leq -(-1+b_1'+b_2')$ which gives $d'-1\leq
-(b_1'+b_2')$, we obtain:
\begin{eqnarray*}
\begin{array}{l}
C^2=(\frac{1}{2}(L-N_1-N_2)+b_1'N_1+b_2'N_2)^2=(d'-1)-2(b_1'^2+b_2'^2)+2(b_1'+b_2')\\
\leq (d'-1)-2(b_1'^2+b_2'^2)-2(d'-1),\\
\end{array}
\end{eqnarray*}
so we obtain
\begin{eqnarray*}
-(d'-1)-2(b_1'^2+b_2'^2)< -2,
\end{eqnarray*}
since $b_1'$, $b_2'$ are integers not both zero. This shows that $C^2=-2$ is not possible in this case.\\
{\bf 2. $d'=1$}. The inequality \eqref{curve-2} becomes
$C^2\leq-2\sum_{i=3}^{8}b_i^2$. If $C^2=-2$ then $-\sum_{i=3}^8 b_i^2=-1$  (so
$0\leq |b_i|\leq 1$), and all the inequalities in \eqref{curve-2} are
equalities. This case corresponds to $C^2=-2$. If there is one index $i$ s.t.
$b_i=-1/2$ then all the $b_i$ are $-1/2$ (because of the structure of the
N\'eron Severi lattice, cf. Corollary \ref{corollary: two divisible classes}),
but then
\begin{eqnarray*}
 -\sum_{i=3}^8 b_i^2=-3/2\not= -1,
\end{eqnarray*}
so $C^2\not= -2$. If one of the $b_i$, $i=3,\ldots,8$ is equal to
$-1$, we have $C=(L-N_1-N_2)/2-N_i$ with $C^2=-2$ and
$C\cdot((L-N_1-N_2)/2)=0$, in fact the divisor has self intersection
zero. In conclusion, for each effective $(-2)$-curve  with
$C\not=N_3,\ldots ,N_8$, we have shown that
\begin{eqnarray*}
C\cdot((L-N_1-N_2)/2)\geq 0.
\end{eqnarray*}
The equality can hold only if $d'=1$. In this case $L^2=4$,
$((L-N_1-N_2)/2)^2=0$ and the divisor is nef but not ample.\\
We show that if the divisor $\bar{D}$ is pseudo ample then $|\bar{D}|$ does not
have base components. As in the previous case
if $|\bar{D}|$ has base components then $\bar{D}=aE+\Gamma$ and we use the same notations as before.\\
{\bf a)} $\Gamma\not= N_i$ for each $i$. Then $0<\Gamma \bar{D}=a-2$ so $a>2$. We have also
$1=\bar{D}N_i=aEN_i+\Gamma N_i$, for $i=1,\ldots,r$ and so
 $E N_i=0$, $\Gamma N_i=1$, for $i=1,\ldots,r$. We have also $EN_j=0$ for each $j=r+1,\ldots,8$. But this is not possible as in the case {\bf a)} above.\\
{\bf b)} $\Gamma=N_i$. {\bf $b_1$)} $\Gamma=N_i$, $i=1,\ldots,r$. We may assume
$i=1$, and so $1=\bar{D}N_1=aEN_1-2$ so $a=3$. Now $1=\bar{D} N_2=3EN_2$ which
is not possible. {\bf $b_2$)} $\Gamma=N_{r+1}$, then $0=\bar{D}N_{r+1}=a-2$, so
$a=2$. On the other hand $1=\bar{D}N_1=2 E N_1$ which is not possible.\\
If $\bar{D}$ is nef and $\bar{D}^2=0$ then by \cite[Theorem
p.79]{kollar} $\bar{D}=kE$, where $|E|$ is a free pencil with
$E^2=0$. By the structure of the N\'eron Severi group we have
$k=1$. \eprf
\begin{cor}
Let $D$ and $\bar{D}$ be divisors as in the Proposition
\ref{prop:compobase}. We suppose that $D^2>0$, $\bar{D}^2>0$.
Let $C$ be a $(-2)$-curve with $C\cdot D=0$ or $C\cdot \bar{D}=0$.
Then $C=N_i$ for some $i=1,\ldots,8$.
\end{cor}

\begin{lemma}\label{lemma:birational2}
With the same notation as in Proposition \ref{prop:compobase}, we have:
\begin{itemize}
\item for $NS(X)=\cl_{2d}$ and $D^2\geq 4$ the map
$\phi_{D}$ is birational, \item for $\bar{D}^2\geq 4$ the
map $\phi_{\bar{D}}$ is birational.
\end{itemize}
\end{lemma}
\bprf The proof is very similar to the proof of Lemma
\ref{lemma:modelli} and is left to the reader.\eprf

We prove the following proposition which is a generalization of
\cite[Proposition 2.6]{catanese} to the case of surfaces in
$\PP^n$.
\begin{prop}\label{prop:cata}
Let  $F$ be a surface in $\PP^n$ and let $\mathcal{N}$ be a subset
of the set of nodes of $F$. Let $G\subset\PP^n$ be an hypersurface
s.t. $div_F(G)=2C$ (here $div_F(G)$ denotes the divisor cut out by $G$ on $F$),
with $C$ a divisor on $F$ which is not Cartier
at the points of $\mathcal{N}$. Then $\mathcal{N}$ is an even set
of nodes iff $G$ has even degree. Conversely if $\mathcal{N}$ is
an even set of nodes then there is an hypersurface $G$ as above.
\end{prop}
\bprf
The proof is identical to the proof of \cite[Proposition 2.6]{catanese}, we recall it briefly.\\
Let $\tilde{F}\longrightarrow F$ be the minimal resolution of the
singularities of $F$, let $H$ denote the pull-back on $\tilde{F}$ of the
hyperplane section on $F$ and let $\deg G=m$ then
\begin{eqnarray*}
mH\sim 2\tilde{C}+\sum\alpha_i N_i
\end{eqnarray*}
where $\tilde{C}$ denotes the strict transform of $C$ on
$\tilde{F}$ and the $N_i$'s denote the exceptional curves over the nodes. Since
$C$ is not Cartier at the singular points, the $\alpha_i$ are odd.
Hence
\begin{eqnarray*}
\sum N_i\sim \delta
H+2\left([\frac{m}{2}]H-\tilde{C}-\sum[\frac{\alpha_i}{2}]N_i\right)
\end{eqnarray*}
where $\delta=0,1$ according to $m$ even or odd. Now if $\sum N_i$
is an even set then $\delta=0$ and $m$ is even. If $m$ is even
then $\delta=0$ and so $\sum N_i$ is an even set.\\
On the other hand, if $\mathcal{N}$ is an even set then
\begin{eqnarray*}
2B\sim\sum N_i.
\end{eqnarray*}
For $B\in Pic(\tilde{F})$ choose $r$ such that $rH-B$ is linearly
equivalent to an effective divisor $\tilde{C}$. Then
\begin{eqnarray*}
2rH\sim 2B+2\tilde{C}=2\tilde{C}+\sum N_i
\end{eqnarray*}
so there is a hypersurface $G$ ($\sim 2rH$) with the properties of
the statement. \eprf From this follows a geometrical
characterization of even set of nodes on K3 surfaces.
\begin{cor}\label{corollary: even set and quadrics}
Let $\bar{X}\subset\PP^{d+1}$, $d\geq 2$, be a surface of degree
$2d$ with a set of $m$ nodes $\mathcal{N}$, s.t. its minimal
resolution is a K3 surface. Then\\
(i)  if $m=8$, then $\mathcal{N}$ is even iff $G$ is a quadric,\\
(ii) if $m=16$ and $d\geq 3$, then $\mathcal{N}$ is even iff $G$
is a quadric.
\end{cor}
\bprf $(i)$ Let $L:=H\cap \bar{X}$ be the generic hyperplane
section with $2d=L^2$ and let $\mathcal{N}$ be an even set of
nodes. Then the lattice $\Z L\oplus N\subset NS(X)$ and we have
$2\hat{N}\equiv\sum N_i$. Since the self-intersection of
$L-\hat{N}$ is $2d-4\geq 0$ by the theorem of Riemann-Roch
$L-\hat{N}$ or $-(L-\hat{N})$ is effective. Since
$(L-\hat{N})\cdot L\geq 0$, $L-\hat{N}$ is effective, so $2L\equiv
2(L-\hat{N})+\sum N_i$. And so $G\in|2L|$ and
$div_{\bar{X}}(G)=2C=2(L-\hat{N})$. The converse follows from the
Proposition \ref{prop:cata}.\\
$(ii)$ The proof in this case is essentially the same. We use the
Kummer lattice $K$ instead of the Nikulin lattice $N$. The lattice
$K$ is generated over $\mathbb{Q}$ by the sixteen disjoint
rational $(-2)$-curves $K_1,\ldots,K_{16}$ and it contains the class
$(K_1+\ldots+K_{16})/2$, which we use in the proof above instead
of the class $\hat{N}$ (for a precise definition of the lattice
$K$ see \cite{nik}).\eprf

\textbf{Remark.} In particular this means that if $\bar{X}$ is a
K3 surface with an even set of nodes, then there exists a quadric
cutting a curve on $\bar{X}$ with multiplicity two, passing
through the even set of nodes (this is the condition $div_F(G)=2C$
of the theorem).


\section{Projective models}\label{section:projective models}

In this section we determine projective models of K3
 surfaces with an even set of nodes and Picard number nine. These were already partially
studied by Barth in \cite{barth1}. Here we recover, with different
methods, some of these examples and we discuss many new examples. Observe that some of
the cases that Barth describes require Picard number at least ten (these
are case five and case four in his list (cf. Paragraph \ref{cl8} below)).\\
In Section \ref{section: ampleness} we proved that the divisors $L$,
$L-\hat{N}$, and $(L-N_1-\ldots-N_m)/2$, $m=2$ or $m=4$ or $m=6$ on $X$ define
regular maps. We use these divisors to give projective models of a K3 surface
$X$ with N\'eron Severi group isometric to $\mathcal{L}_{2d}$ or to
$\mathcal{L}_{2d}'$ and in general we study the projective models of the same
surface by using different polarizations. In particular in each case one can
use as polarization $L$ or $L-\hat{N}$, if $L^2>4$. The first polarization
contracts the curves of the even sets to eight nodes on the surface, the second
one sends these curves to lines on the projective model.\\
In case (1) of the Corollary \ref{corollary: NS of X and of Y} it
is also possible to study the projective models given by the maps
$\phi_{L_1}$, resp. $\phi_{L_2}$ (for $L_i^2>0$) where $L_1$ and
$L_2$ are the divisors defined in \eqref{formula: L1 and L2: first
case } or in \eqref{formula: L1 and L2: second case }. They give
projective models of $X$ in the projective space
$\PP(H^0(X,L_1))$, resp. $\PP(H^0(X,L_2))$ or give $2:1$ maps to
the images of $X$ in these spaces. If the maps are not $2:1$, the
image of $X$ contains nodes and lines, which on $X$ form an even
set. The image of $X$ under
$\phi_{L_1}\times\phi_{L_2}:X\rightarrow\PP(H^0(X,L_1))\times
\PP(H^0(X,L_2))$ is a surface, this surface is the image of
$Y\hookrightarrow\PP^{h^0(L_1)+h^0(L_2)-1}$ under the projection
to the eigenspaces: $\PP^{h^0(L_1)+h^0(L_2)-1}\longrightarrow
\PP(H^0(X,L_1))\times\PP(H^0(X,L_2))$. Indeed put
$h^0(L_1)=m_1+1$, $h^0(L_2)=m_2+1$, then $h^0(M)=m_1+m_2+2$ and we
have a commutative diagram:
$$
\begin{xy}
\xymatrix{
Y\ar[r]\ar@{-->}[d]        &\PP^{m_1+m_2+1}           \ar[r]^v        &**[r]\ar@{->}[d]^p\PP^r\supset V_{m_1+m_2+1}\\
X\ar[r]                    &\PP^{m_1}\times\PP^{m_2}  \ar[r]^s
&**[r]\PP^{r'}\supset S_{m_1+m_2} }
\end{xy}
$$
Here the rational map between $Y$ and $X$ follows from Diagram
\ref{diagrammone}, $v$ is the Veronese embedding and so
$r=\frac{(m_1+m_2+2)(m_1+m_2+3)}{2}$, $s$ is the Segre embedding
and so $r'=(m_1+1)(m_2+1)-1$, $p$ is the projection,
$V_{m_1+m_2+1}$ is the image of $\PP^{m_1+m_2+1}$ in $\PP^r$ and
$S_{m_1+m_2}$ is the image of $\PP^{m_1}\times\PP^{m_2}$ in
$\PP^{r'}$. The Nikulin involution on $\PP^{m_1+m_2+1}$  operates
as:
$$
(x_0:\ldots:x_{m_1}:y_0:\ldots:y_{m_2})\mapsto (x_0:\ldots:x_{m_1}:-y_0:\ldots:-y_{m_2})
$$
which induces an operation on the coordinates of $\PP^r$ as
$$
\begin{array}{l}
(x_0^2:x_1^2:\ldots:y_0^2:y_1^2:\ldots:y_{m_2-1}y_{m_2}:x_0y_1:x_0y_2:\ldots:x_{m_1}y_{m_2})\mapsto\\
(x_0^2:x_1^2:\ldots:y_0^2:y_1^2:\ldots:y_{m_2-1}y_{m_2}:-x_0y_1:-x_0y_2:\ldots:-x_{m_1}y_{m_2}).
\end{array}
$$
The projection $p$ goes to the invariant space $\PP^{r'}$ with coordinates
$(x_0y_1:x_0y_2:\ldots:x_{m_1}y_{m_2})$ and so the equations of the image of
$Y$ in  $V_{m_1+m_2+1}$ in these coordinates give equations for the image of
the
surface $X$ in  $S_{m_1+m_2}$ (cf. also \cite[Proposition 2.7]{bertio}).\\
Observe that the sum of the divisors $L_1$ and $L_2$ is exactly
$L-\hat{N}$. From now on if $NS(X)=\cl_{2d}'$ and $d/2$ is odd
then $L_1:=(L-N_1-N_2)/2$, $L_2:=(L-N_3-\ldots-N_8)/2$, if
$NS(X)=\cl_{2d}'$ and $d/2$ is even then
$L_1:=(L-N_1-\ldots-N_4)/2$, $L_2:=(L-N_5-\ldots-N_8)/2$.\\

In the case $NS(X)=\cl_{2d}$ the construction above holds if
instead of $L_1$ and $L_2$ we take $L$ and $L-\hat{N}$.

\begin{footnotesize}
\begin{table}\caption{N\'eron Severi lattices and projective models}\label{tab1}
$$
\begin{array}{lclc}
&X&&Y\\
&NS(X)=\cl_{2}&&NS(Y)=\mathcal{M}_{4}'\\
\phi_L &\mbox{double plane (singular sextic)}&
\phi_M\ &\mbox{smooth quartic in }\mathbb{P}^3\\
\\
&NS(X)=\cl_{4}&&NS(Y)=\mathcal{M}_{8}'\\
\phi_L &\mbox{quartic with even set of nodes}&
\phi_M &\mbox{complete intersection in }\mathbb{P}^5\\
\\
&NS(X)=\cl_{4}'& &NS(Y)=\mathcal{M}_{2}\\
\phi_L &\mbox{double cover of a cone}&\phi_M &\mbox{double plane}\\
\phi_{L_1}&\mbox{elliptic fibration}&&\\

\\
&NS(X)=\cl_{6}& &NS(Y)=\mathcal{M}_{12}'\\
\phi_L\ &\mbox{singular complete intersection in }\mathbb{P}^4&\phi_M&\mbox{projective model in }\mathbb{P}^7\\
\phi_{L-\hat{N}}&\mbox{double plane (smooth sextic)}&&\\
\phi_L\times\phi_{L-\hat{N}}&\mbox{complete intersection in }
\mathbb{P}^4\times\mathbb{P}^2&&\\

\\
&NS(X)=\cl_{8}&&NS(Y)=\mathcal{M}_{16}'\\
\phi_L\ &\mbox{singular complete intersection in }\mathbb{P}^5&\phi_M&\mbox{projective model in }\mathbb{P}^9\\
\phi_{L-\hat{N}}&\mbox{smooth quartic in }\mathbb{P}^3&&\\
\\
&NS(X)=\cl_{8}'& &NS(Y)=\mathcal{M}_{4}\\
\phi_L\ &\mbox{singular complete intersection in }\mathbb{P}^5&\phi_M&\mbox{smooth quartic in }\mathbb{P}^3\\
\phi_{L-\hat{N}}&\mbox{double cover of a quadric}\\

\\
&NS(X)=\cl_{10}&&NS(Y)=\mathcal{M}_{20}'\\
\phi_{L-\hat{N}}&\mbox{smooth complete intersection in }\mathbb{P}^4&\phi_M&\mbox{projective model in }\mathbb{P}^{11}\\
\phi_{L-\sum_{i=1}^4N_i}&\mbox{double cover of a plane}&&\\

\\
&NS(X)=\cl_{12}& &NS(Y)=\mathcal{M}_{24}'\\
\phi_{L-\hat{N}}&\mbox{smooth complete intersection in }\mathbb{P}^5&\phi_M&\mbox{projective model in }\mathbb{P}^{13}\\
\phi_{L-\sum_{i=1}^4N_i}&\mbox{singular quartic in }\mathbb{P}^3&&\\
&\mbox{(mixed even set with conics)}\\
\\
&NS(X)=\cl_{12}'& &NS(Y)=\mathcal{M}_{6}\\
\phi_{L-\hat{N}}&\mbox{smooth complete intersection in }\mathbb{P}^5&\phi_M&\mbox{complete intersection in }\mathbb{P}^{4}\\
\phi_{L_2}\times\phi_{L_1}&\mbox{surface of bidegree } (2,3) \mbox{ in }\mathbb{P}^1\times\mathbb{P}^2&&\\
\\
&NS(X)=\cl_{16}'& &NS(Y)=\mathcal{M}_{8}\\
\phi_{L_1}\times\phi_{L_2}&\mbox{complete intersection in
}\mathbb{P}^2\times\mathbb{P}^2&
\phi_M&\mbox{complete intersection in }\mathbb{P}^{5}\\
\\
&NS(X)=\cl_{24}'& &NS(Y)=\mathcal{M}_{12}\\
\phi_{L_1}\times\phi_{L_2}&\mbox{complete intersection in
}\mathbb{P}^3\times\mathbb{P}^3&
\phi_M&\mbox{complete intersection in }\mathbb{P}^{7}
\end{array}
$$
\end{table}
\end{footnotesize}

\subsection{The case of $L^2=2$, $NS(X)=\cl_2$, the polarization
$L$}\label{cl2} Since $L$ is pseudo ample by the Proposition
\ref{prop:pseudoample} the linear system $|L|$ defines a $2:1$ map $X'\lra
\PP^2$ ramified on a sextic curve with eight nodes where $X'$ is the surface
$X$ after contraction of the $(-2)$-curves. More precisely we have a
commutative diagram:
\begin{eqnarray*}
\begin{array}{ccc}
X&\lra&X'\\
\downarrow&&\downarrow\\
\widetilde{\PP}^2&\lra&\PP^2\\
\end{array}
\end{eqnarray*}
where $\widetilde{\PP}^2$ is the blow up of $\PP^2$ at the eight
double points of the sextic. By general results 
on cyclic coverings the pull back of the branching sextic on $X$ is
$3L-(N_1+\ldots+N_8)=3L-2\hat{N}=(L-\hat{N})+(2L-\hat{N})$. Now
$(L-\hat{N})^2=2-4=-2$, $L\cdot (L-\hat{N})=2$  and so by using Riemann-Roch
Theorem the divisor $L-\hat{N}$ is effective, and is a rational curve of degree
two on $\PP^2$. Observe that its image is an irreducible conic, in fact we are
assuming that $X'$ has exactly eight nodes and no other singularities. On the
other hand $(2L-\hat{N})^2=8-4=4$ so by Proposition \ref{prop:compo1} and  by
\cite[Proposition 2.6]{saintdonat} the generic member in $|2L-\hat{N}|$ is an
irreducible curve of genus three, and its image in $\PP^2$ is a curve of degree
four (and in fact genus $3=(4-1)(4-2)/2$).  In both cases we have
$(L-\hat{N})\cdot N_i=(2L-\hat{N})\cdot N_i=1$ and so the curves intersect at
the points which are the images of the curves $N_i$ in $\PP^2$. This is the
first case in the paper of Barth, \cite{barth2}.
\subsection{The case of $L^2=4$, $NS(X)=\cl_4$, the polarization $L$}\label{cl4} By
the Proposition \ref{prop:pseudoample} the linear system $|L|$ defines a
birational map $\phi_L$ from $X$ to a quartic surface in $\PP^3$, the curves
$N_i$ are contracted to nodes. In this case $(L-\hat{N})^2=0$ and by
 the Proposition \ref{prop:compo1} $|L-\hat{N}|$ has no base components, so the generic member in the system is
an irreducible elliptic curve (observe that by the structure of
the N\'eron Severi group it cannot be the multiple of an elliptic
curve). Since $L\cdot (L-\hat{N})=4$ the elliptic curve is sent to
a quartic curve in $\PP^3$  and is a complete intersection of two
quadrics (observe that it cannot be a plane quartic since this has
genus three). Moreover since $(L-\hat{N})\cdot N_i=1$, the quartic
contains the nodes. There is a third quadric passing through the
nodes, in fact $h^1(L-\hat{N})=0$ (\cite[Proposition
2.6]{saintdonat}) hence $h^0(L-\hat{N})=2$ and again by loc. cit.
$h^0(2(L-\hat{N}))=3$. Now $2(L-\hat{N})=2L-(N_1+\ldots+ N_8)$ and
the image of these divisors are precisely the quadrics which
vanish on the eight singular points (cf. Corollary \ref{corollary:
even set and quadrics}). Let $s_1,s_2$ be a basis of
$H^0(L-\hat{N})$ then $s_1^2,s_1s_2,s_2^2$ is a basis of
$H^0(2(L-\hat{N}))$ and these are the three quadrics through the
nodes. This is the case three of Barth \cite{barth2} and by the
Table \ref{tab1} it corresponds to the case of $\mathcal{M'}_{8}$
of \cite{bertio}.
\subsection{The case of $L^2=4$, $NS(X)=\cl'_4$} \label{cl4'}
\begin{itemize}
\item[{\bf(a)}] {\bf The polarization $L$}. We may assume that the
class $(L/2,v/2)$ is equal to $(L/2,(-N_1-N_2)/2)$. By the
Proposition \ref{prop:compobase} and \cite[Proposition
2.6]{saintdonat} this defines a pencil of elliptic curves which we
denote by $E$. Observe that $L=2E+N_1+N_2$ with $N_i\cdot E=1$,
hence by \cite[Proposition 5.7, (iii), a)]{saintdonat} $L$ defines
a $2:1$ map to a cone of $\PP^3$. The pencil $|E|$ corresponds to
the system of lines through the vertex of the cone under this map.
The class $C_2:=E-\hat{N}+N_1+N_2=L/2+(-N_3-\ldots-N_8)/2$ is
effective with $C_2\cdot N_i=1$, similarly  we may assume that the
class $C_6:=3L/2+(-N_3-\ldots-N_8)/2$ is an irreducible curve (it
follows by Proposition \ref{prop:compobase}), with $C_6\cdot
N_i=1$, moreover $C_2\cdot C_6=0$, $C_2\cdot L=2$ and $C_6\cdot
L=6$. Let $c_2:=\varphi_L(C_2)$ and $c_6:=\varphi_L(C_6)$. These
two curves meet on the cone at the images of $N_i$, $i=3,\ldots,
8$. Their union is a curve of degree eight, which is the branch
divisor of the covering. In fact if $C_2$ is not a component of
the branch divisor then $\varphi_L(C_2)$ has degree one and so is
a line. But this means that $C_2\in|E|$  which is not the case.
Hence $C_2$ is a component of the branch divisor and $c_2$ is a
conic. If now $C_6$ is not in the branch divisor, we have $\deg
c_6=3$, and $c_2\cdot c_6=6$, but then $c_6$ is contained in the
plane of $c_2$ and on the cone too, which is impossible. Hence
$C_6$ is also a component of the branch divisor. Finally observe
that $\varphi_L(N_i)=Q$, for $i=1,2$ where $Q$ is the vertex of
the cone. This surface is also described in \cite[Paragraph
3.2]{bertio}. \item[{\bf (b)}] {\bf An elliptic fibration}. Now we
describe the elliptic fibration on $X$ defined by the divisor $E$.
We consider the rational curve $C_2=L/2+(-N_3-\ldots-N_8)/2$,
it has intersection one with the class of the fiber $E$. So $C_2$
is a section of the fibration $\phi_E:X\rightarrow \mathbb{P}^1$
and
the classes $E$ and $C_2$ generate a lattice isometric to $U$.\\
Since the six $(-2)$-curves $N_3,\ N_4,\ldots,N_8$ are orthogonal
to $E$, they are the components of some reducible fibers. All
these curves intersect the section $C_2$ so they are components of
six different reducible fibers. The rational curve $N_1$ is
another section of the fibration (because its intersection with
$E$ is one). The N\'eron Severi group is generated over $\Q$ by
the classes $E$ of the fiber, by $C_2$, by the components $N_i$,
$i=3,\ldots,8$ of the reducible fibers and by the other section
$N_1$. The N\'eron Severi group of an elliptic fibration admitting
a section is generated by the class of the fiber, by the zero
section, by the irreducible components of the reducible fibers
(not meeting the zero section) and by other sections. Since the
Picard number is nine the six reducible fibers containing $N_i$,
$i=3,\ldots, 8$ are the only reducible fibers of the fibration,
they are all of type $I_2$ (two rational curves meeting in two
distinct points). The Euler characteristic of a K3 surface is 24
and is the sum of the Euler characteristics of the singular
fibers. The singular irreducible fibers in the generic case are of
type $I_1$ (singular irreducible curve with a node). Each fiber of
type $I_1$ has Euler characteristic one, and each fiber of type
$I_2$ has Euler characteristic equal to two. By the computation on
the Euler characteristic it is clear that there are twelve
singular fibers of type $I_1$ and six of type $I_2$. There are two
independent sections, so the rank of the Mordell-Weil lattice is
one. One of these sections (the zero section) is the curve $C_2$,
mapped by $\phi_{L}$ to the conic in the branch locus on the cone.
Other sections correspond to the curves
$N_1$ and $N_2$, these are both mapped to the vertex of the cone.\\
Since $X$ is a double cover of a cone, it admits an involution
$j$. This involution fixes the classes $N_i,$ $i=3,\ldots, 8$,
because they correspond to the intersection points between the
conic and the sextic in the branch locus on the cone; it fixes the
class $C_2$, and switches the classes $N_1$ and $N_2$. The
involution $j$ fixes also the class $L$ which defines the double
cover $\phi_L$. On the fibration the involution $j$ fixes the
class of the fiber $E$ and so it acts on the base, $\PP^1$, of the
fibration as the indentity, it fixes the zero section, which
corresponds to the class $C_2$, and switches the other two
independent sections $N_1$ and $N_2$. On the reducible fiber the
involution $j$ clearly fixes the component $N_i$ and it fixes the
other component $E-N_i$ too, since it fixes the fiber and a
reducible fiber has two components. On $E-N_i$ the involution $j$
switches the points $P_1$ and $P_2$, which are the points of
intersection between the fiber
and the sections $N_1$, respectively $N_2$.\\
The surface $X$ admits an even set of eight disjoint rational
curves, so it is the minimal resolution of the quotient of a K3
surface $Y$ by a Nikulin involution. The elliptic fibration of $X$
on $\mathbb{P}^1$ induces a fibration of $\widetilde{Y}$ (the blow
up of $Y$) on $\mathbb{P}^1$ and so of $Y$ (cf. Diagram
\ref{diagrammone}). Let $E$ denote the generic fiber of the
fibration on $X$ and $A$ the generic fiber of the fibration on
$\widetilde{Y}$. By the Hurwitz formula, we have
$$2g(A)-2=2(2g(E)-2)+\deg R$$
where $R$ is the branch divisor. Since $X$ has an elliptic fibration we have
$g(E)=1$ and $\deg R=2$ because the involution ramifies on the points of
intersection $E\cap N_1$ and $E\cap N_2$. So we find $2g(A)-2=2(2-2)+2$, hence
the generic fiber of the fibration $\widetilde{Y}\rightarrow\mathbb{P}^1$ is
hyperelliptic of genus two.
\end{itemize}
\subsection{The case of $L^2=6$, $NS(X)=\cl_6$}\label{cl6}
\begin{itemize}
\item[{\bf (a)}] {\bf The polarization $L-\hat{N}$}. In this case
$(L-\hat{N})^2=2$ by Lemma \ref{lemma:modelli} it defines a $2:1$ map to
$\PP^2$. The curves $N_i$ are mapped to lines in the plane. Let
$l:=\phi_{L-\hat{N}}(L-\hat{N})$, then for each curve $C$ in the plane we have
the formula $\phi_{L-\hat{N}}^*(l)\cdot\phi_{L-\hat{N}}^*(C)=2(C\cdot l)$.
Since $(L-\hat{N})\cdot N_i=1$ the curves $N_i$ are contained in the preimage
$\phi_{L-\hat{N}}^{-1}(T_i)$ where $T_i$ are lines which are tritangents to the
branch divisor, and so $\varphi_{L-\hat{N}}^*(T_i)=N_i+N_i'$. The curves $N_i,
N_i'$
 meet in three points. Barth in \cite[Paragraph 2]{barth2}
shows that there is a quartic in $\PP^2$ meeting the branch
sextic at the tangency points.\\ 
\item[{\bf (b)}] {\bf The polarization $L$}. We consider the projective model
of $X$  as complete intersection of a cubic and a quadric hypersurface, it has
eight nodes and the map is $\phi_{L}:X\rightarrow\mathbb{P}^4$. The curve
$L-\hat{N}$ (cf. Proposition \ref{prop:compo1}) has degree
$6=L\cdot(L-\hat{N})$ and genus $(L-\hat{N})^2/2+1=2$. Since $(L-\hat{N})\cdot
N_i=1$, $i=1,\ldots,8$ its image in $\mathbb{P}^4$ passes through the eight
singular points. This curve is contained in the intersection of three quadrics
in $\mathbb{P}^4$, in fact $h^0(2L-(L-\hat{N}))=(L+\hat{N})^2/2+2=3$. The eight
singular points of the surface are contained in three more quadrics, in fact
$h^0(2L-(\sum_{i=1}^8N_i))=6$ (cf. Corollary \ref{corollary: even set and quadrics}).\\
We consider now the linear system $|L-\hat{N}|$ associated to the
hyperplane sections passing through the eight
singular points of the image of $X$ in $\mathbb{P}^4$. We have $h^0(L-\hat{N})=3$ and
let $l_1$, $l_2$, $l_3$ be its generators. The six elements
$l_1^2$, $l_2^2$, $l_3^2$, $l_1\cdot l_2$, $l_1\cdot l_3$,
$l_2\cdot l_3$ span $|2L-\sum_{i=1}^8N_i|\cong \PP^5$ (these are the quadrics passing through the nodes).\\
\item[{\bf (c)}] {\bf The map $\phi_L\times\phi_{L-\hat{N}}$}. In
\cite[Paragraph 3.9]{bertio} the K3 surface $Y$ admitting a Nikulin involution
with N\'eron Severi group $\mathcal{M}_{12}'$ is described. Its quotient
$\bar{X}$ is birational to a K3 surface which is complete intersection of a
hypersurface of bidegree $(2,0)$ and three hypersurfaces of bidegree $(1,1)$ in
$\mathbb{P}^4\times \mathbb{P}^2$ (for a more detailed description of this
complete intersection see the Section \ref{geometric}). The minimal resolution
of the quotient $\bar{X}$ is the K3 surface $X$ with $NS(X)=\mathcal{L}_{6}$.
The projection of $X$ to the first factor is defined by the divisor $L$ and to
the second one by $L-\hat{N}$. The first projection contracts eight disjoint
rational curves and the same curves are sent to eight lines by the second
projection (the $2:1$ map to $\mathbb{P}^2$).
\end{itemize}

\subsection{The case of $L^2=8$, $NS(X)=\cl_8$}\label{cl8}
\begin{itemize}
\item[{\bf (a)}] {\bf The polarization $L-\hat{N}$}. We have $(L-\hat{N})^2=4$
and the map $\phi_{L-\hat{N}}:X\lra\PP^3$ exhibits $X$ as a quartic surface in
$\PP^3$ with eight disjoint lines. This case is studied by Barth in
\cite{barth2}. He describes two conditions to have an even set. The second one
is not satisfied in our case, since it requires Picard number at least ten. In
fact he shows that in this case there are two skew lines $Z_1$, $Z_2$ on the
quartic surface with $Z_1$ meeting four lines and skipping the other four
lines, and viceversa for $Z_2$. An easy computation shows that the intersection
matrix of
the hyperplane section, of the lines $N_i$ and of $Z_1$ (or $Z_2$) has rank ten.\\
Barth's first condition
says that there is an elliptic quartic curve in $\PP^3$ which meets in two points
four rational curves and skips the other four. In term of classes in the N\'eron Severi
lattice this means that there is a curve $E=\alpha L+\sum b_i N_i$ with $E^2=0$, $E N_i=2$
 for $i=1,\ldots,4$ and $EN_i=0$ for $i=5,\ldots,8$. By using the intersection products we obtain
 $\alpha=1$, $b_i=-1$ for $i=1,\ldots,4$ and $b_i=0$ for $i=5,\ldots,8$. So the elliptic curve is
  $E=L-N_1-\ldots-N_4$ and in fact $E\cdot (L-\hat{N})=4$ which is the degree of $E$ in $\PP^3$.
Similarly the curve $L-N_5-\ldots-N_8$ meets the other four curves and skips the first four. Finally observe that these divisors are not studied in Proposition \ref{prop:compo2}, with the notation there this is the case $d=r$ and the proof does not work in this case.\\
We describe briefly the elliptic fibration defined by $E$. Since $E\cdot
N_i=0$, for $i=5,\ldots,8$ these are components of reducible fibers. On the
other hand the curves $N_i$, $i=1,\ldots,4$ are bisections of the fibration.
The curves $L-N_2-N_3-N_4-N_j-N_k$ with $j\not=k$ and $j,k=5,6,7,8$ are
rational $(-2)$-curves which meet $E$ in two points and $N_j$, $N_k$ in two
points as well. Hence they are also bisections of the fibration, and since a
bisection meets also the singular fiber in two points, the curves
$N_5,\ldots,N_8$ are contained in four different singular fibers, which are of
type $I_2$. The remaining singular fibers are of type $I_1$, and we have $16$
of them. This fibration does not admit sections. In this case the even set
consists of four bisections and of four components of the singular fibers $I_2$
(these are all disjoint). \item[{\bf (b)}] {\bf The polarization $L$}. We
consider the projective model of $X$ given by the map\\
$\phi_{L}:X\rightarrow\mathbb{P}^5$, this is a complete
intersection of three quadrics and has eight nodes. The generic
element in $|L-\hat{N}|$ is a curve of degree
$8=L\cdot(L-\hat{N})$ and genus $(L-\hat{N})^2/2+1=3$ (cf.
Proposition \ref{prop:compo1}). Since $(L-\hat{N})\cdot N_i=1$,
$i=1,\ldots,8$ the image of the curve $L-\hat{N}$ in
$\mathbb{P}^5$ passes through the eight singular points, moreover
this divisor is not Cartier at the nodes. By the Corollary
\ref{corollary: even set and quadrics} there exists a quadric $G$
which cuts on the surface the curve $L-\hat{N}$ passing through
all the singular points, so this curve is contained in the
intersection of four quadrics in $\mathbb{P}^5$, in fact
$h^0(2L-(L-\hat{N}))=(L+\hat{N})^2/2+2=4$. The quadric $G$ must
cut the image of $L-\hat{N}$ with multiplicity two since
$\deg(X)=8$ and the intersection has degree $16$. The eight
singular points of the surface are contained in the intersection
of ten quadrics, in fact $h^0(2L-(\sum_{i=1}^8N_i))=10$. We
consider the linear system $|L-\hat{N}|$ associated to the
hyperplane section passing through the eight singular points. We
have $h^0(L-\hat{N})=4$ and we call $l_1$, $l_2$, $l_3$, $l_4$ its
generators. The ten elements $l_1^2$, $l_2^2$, $l_3^2$, $l_4^2$,
$l_1l_2$, $l_1l_3$, $l_1l_4$, $l_2l_3$, $l_2l_4$, $l_3l_4$ span
$|2L-\sum_{i=1}^8N_i|$.
\item[{\bf (c)}] {\bf The map $\phi_L\times\phi_{L-\hat{N}}$}. 
This K3 surface is the minimal resolution of the quotient of a K3
surface $Y$ by a Nikulin involution. The N\'eron Severi group of
$Y$ is $\mathcal{M}_{16}'$ by the Table \ref{tab1}, and $M$ is the
ample class on $Y$ with $M^2=16$. This gives an immersion of $Y$
in $\PP^9$, and the action of the Nikulin involution is induced by
$(x_0:\ldots:x_5:y_0:\ldots:y_3)\mapsto(x_0:\ldots:x_5:-y_0:\ldots:-y_3)
$. By the projection formula we have
\begin{eqnarray*}
H^0(Y,M)\cong H^0(X,L)\oplus H^0(X,L-\hat{N}),
\end{eqnarray*}
with $h^0(X,L)=6$, $h^0(X,L-\hat{N})=4$. Now
\begin{eqnarray*}
S^2H^0(Y,M)=(S^2H^0(X,L)\oplus S^2H^0(X,L-\hat{N}))\oplus(H^0(X,L)\otimes
H^0(X,L-\hat{N})).
\end{eqnarray*}
This has dimension $55=(21+10)+24$. On the other hand
\begin{eqnarray*}
H^0(Y, 2M)\cong H^0(X,2L)\oplus H^0(X,2L-\hat{N})
\end{eqnarray*}
and the dimensions are $34=18+16$. This shows that there are $(21+10)-18=13$
invariant quadrics and $24-16=8$ antiinvariant quadrics $Q_i(x,y)$, $i=1\ldots,
8$ in the ideal of $Y$. Since the quadrics in four variables are only ten,
there are three quadrics
$q_1(x_0,\ldots,x_5),~q_2(x_0,\ldots,x_5),~q_3(x_0,\ldots,x_5)$ in the ideal of
$Y$. The map $\phi_L\times \phi_{L-\hat{N}}$ sends $X$ to the product
$\PP^5\times\PP^3$ and its image is the image of $Y\subset \PP^9$ into the
product of the eigenspaces, hence it is contained in three quadrics
$q_1(x_0,\ldots,x_5),~q_2(x_0,\ldots,x_5),~q_3(x_0,\ldots,x_5)$ of bidegree
$(2,0)$ and eight quadrics $Q_i(x,y)$, $i=1,\ldots, 8$, of bidegree $(1,1)$ in
particular it is not a complete intersection of quadrics. The quadrics
$q_1(x_0,\ldots,x_5),~q_2(x_0,\ldots,x_5),~q_3(x_0,\ldots,x_5)$ define the
image of $Y$ in $\PP^5$, which is $X$ with the polarization $L$. Since the
fixed points of the Nikulin involution are contained in the space
$y_0=\ldots=y_3=0$, then the projection of the ten quadrics of the kind
$q(x)-q'(y)=0$ to $\PP^5$  are ten quadrics cutting out the set of nodes on
$X\subset \PP^5$. The projection to $\PP^3$ is $X$ with the polarization
$L-\hat{N}$ and is a quartic. One can obtain an equation for the quartic in the
following way: a point $x\in X\subset \PP_5$ has a non-trivial counterimage if
there is a non-trivial solution of $Q_i(x,y)=\sum_{j=0}^{5} a_{ij}(y)x_j=0$,
$i=1,\ldots, 8$ which for a fixed $x$ is a linear system of eight equations in
six variables. Hence all the $6\times 6$ minors of the matrix $(a_{ij}(y))$ are
zero. Each of these is a sextic surface of $\PP^3$ vanishing on $X\subset
\PP^3$. Since this is a surface of degree four, each of them splits into a
product $q(x)\cdot p_4(x)$ where $p_4(x)=0$ is an equation of $X\subset \PP^3$.
\end{itemize}

\subsection{The case of $L^2=8$, $NS(X)=\cl_8'$}\label{cl8'}
\begin{itemize}
\item[{\bf (a)}] {\bf The polarization $L-\hat{N}$}. In this case we have the
divisor $E_1:=(L/2, v/2)$ with $v^2=-8$ and $v=-N_1-N_2-N_3-N_4$, and also the
divisor $E_2:=(L/2,v'/2)$ with $v'=-N_5-N_6-N_7-N_8$ so
$(L/2,v/2)^2=(L/2,v'/2)^2=0$ and $L-\hat{N}=(L/2,v/2)+(L/2,-v/2)$ is the sum of
two elliptic curves (cf. Proposition \ref{prop:compobase}). This is a $2:1$ map
to $\PP^1\times \PP^1$ (by Lemma \ref{lemma:modelli}) and the curves $N_i$ are
sent to lines on the quadric. Moreover since $E_1\cdot N_i=1$ and $E_2\cdot
N_i=0$ for $i=1,2,3,4$ the images of these lines belong to the same ruling on
the quadric and the images of $N_i$, $i=5,6,7,8$ belong to the other ruling. By
a similar computation as in Paragraph \ref{cl6} the curves $N_i$ are one of the
two components of the preimage of a curve on the quadric which splits on $X$,
hence $\phi_{L-\hat{N}}(N_i)=T_i$ and these are bitangents to the branch curve
of bidegree $(4,4)$. Let $\phi_{L-\hat{N}}^*(T_i)=N_i+N_i'$ then
\begin{eqnarray*}
\mathcal{O}_{B'}(N_1+\ldots
+N_8+N_1'+\ldots+N_8')=\mathcal{O}_{B'}(4(L-\hat{N}))=\mathcal{O}_{B'}(2(2(L-\hat{N})))
\end{eqnarray*}
By Proposition \ref{prop:compo1} $2(L-\hat{N})$ is a curve, and has bidegree
$(2,2)$ on the quadric. Hence the divisor cut out by $N_1+\ldots+N_8'$ is two
times the divisor cut out by $2(L-\hat{N})+\mu$ where $\mu$ is a 2-torsion
element in the Picard group. Barth shows in \cite{barth2}, case six, that such
an element does not exist. This implies that the tangency points of the $T_i$
on the quadric are cut out by a curve of bidegree $(2,2)$.
\item[{\bf (b)}] {\bf The polarization $L$}. All the considerations of the
Paragraph \ref{cl8}, case {\bf (b)} are true. Moreover there are two elliptic
curves $L_1=\frac{L-N_1-N_2-N_3-N_4}{2}$, $L_2=\frac{L-N_5-N_6-N_7-N_8}{2}$
passing through four of the eight singular points each and not passing through
the other four. Obviously also in this case the image of
$2(L-\hat{N})=2(L_1+L_2)$ is cut out by a quadric. By the Table \ref{tab1} this
case corresponds to a K3 surface $Y$ with $NS(Y)=\mathcal{M}_4$. After a change
of coordinates the surface $X$  can be written in the form
\begin{eqnarray*}
q(z_0,\ldots,z_5)=0,\,\, z_0z_1-z_4^2=0,\,\, z_2z_3-z_5^2=0
\end{eqnarray*}
and there are four singularities on $z_0=z_1=z_4=0$ and four on
$z_2=z_3=z_5=0$ (the two copies of $\PP^2$ which are the vertices
of the cones). Now the quadrics of the kind $z_iz_j=0$ with
$i=0,1$ and $j=2,3$ meet the K3 surface in two curves $C_i$, $C_j$
with multiplicity two, hence $2C_i\in|L-(N_1+\ldots+N_4)|$ and
$2C_j\in|L-(N_5+\ldots+N_8)|$ and so $C_i$, $C_j$  are in the
linear system of $L_1$, resp. of $L_2$.
\end{itemize}

\subsection{The case of $L^2=10$, $NS(X)=\cl_{10}$}\label{cl10}

\begin{itemize}
\item[{\bf (a)}] {\bf The polarization $L-\hat{N}$.} Since $(L-\hat{N})^2=6$,
then the projective model of $X$ is a complete intersection of a quadric and a
cubic hypersurfaces in $\mathbb{P}^4$ with an even set of eight lines (cf.
Lemma \ref{lemma:modelli}). \item [{\bf (b)}]{\bf The polarizations
$L-N_1-N_2-N_3-N_4$ and $L-N_5-N_6-N_7-N_8$}. The divisors $L-N_1-N_2-N_3-N_4$
and $L-N_5-N_6-N_7-N_8$ are pseudo ample classes by the Proposition
\ref{prop:compobase}. They define two maps $2:1$ to $\mathbb{P}^2$. Each of
these maps contracts four curves of the eight rational curves $N_i$ and maps
the other four in four conics.
\end{itemize}

\subsection{The case of $L^2=12$, $NS(X)=\cl_{12}$}\label{cl12}

\begin{itemize}
\item [{\bf (a)}]{\bf The polarization $L-\hat{N}$.} Since $(L-\hat{N})^2=8$
the projective model of $X$ is a K3 surface in $\mathbb{P}^5$ with an even set
of eight disjoint lines. \item [{\bf (b)}]{\bf The polarizations
$L-N_1-N_2-N_3-N_4$ and $L-N_5-N_6-N_7-N_8$}. The curves
$E_1=L-N_1-N_2-N_3-N_4$ and $E_2=L-N_5-N_6-N_7-N_8$ have self intersection
four, so they define two maps to $\mathbb{P}^3$ (by Lemma
\ref{lemma:birational2}). The map $\phi_{E_1}$ contracts the four curves $N_i$,
$i=1,\ldots, 4$ and sends the other in four conics. The map $\phi_{E_2}$
contracts the other four curves and sends $N_i$, $i=1,\ldots,4$  in conics.
\end{itemize}

\subsection{The case of $L^2=12$,
$NS(X)=\cl_{12}'$}\label{cl12'}

\begin{itemize}
\item [{\bf (a)}]{\bf The polarization $L-\hat{N}$.}
Observe that the considerations of Paragraph
\ref{cl12}, {\bf (a)} are true also in this case. Moreover there
are two curves $C_1=\frac{L-N_1-N_2-N_3-N_4-N_5-N_6}{2}$ and
$C_2=\frac{L-N_7-N_8}{2}$ intersecting respectively six and two of
the lines $N_i$ in one point. The curve $C_1$ has degree three and
genus
one. The curve $C_2$ has degree five and genus two.
\item [{\bf (b)}]{\bf The map $\phi_{L_1}\times \phi_{L_2}$}.  The intersection
properties of $L_1$ and $L_2$ are $L_1\cdot L_1=2$, $L_2\cdot
L_2=0$ and $L_1\cdot L_2=3$. The K3 surface $X$ is the minimal rsolution of the quotient $\bar{X}$ of
a K3 surface $Y$ admitting a Nikulin involution with
$NS(Y)=\mathcal{M}_{6}$ which is described in \cite[Paragraph 3.3]{bertio}. The surface $\bar{X}$ has bidegree $(2,3)$ in $\PP^1\times \PP^2$. The maps $\phi_{L_1}$ and $\phi_{L_2}$ are respectively
the projection to the second and to the first projective space.\\
The map $\phi_{L_1}:X\rightarrow\mathbb{P}^2$ is a 2:1 map. It
contracts the six rational curves $N_3,\ldots, N_8$ to six nodes
of the branch sextic and the two curves $N_1$ and $N_2$ are mapped
to lines in $\mathbb{P}^2$ which are tritangent to the branch
locus. The map $\phi_{L_2}:X\rightarrow\mathbb{P}^1$ is an
elliptic fibration, it contracts the two rational curves $N_1,\
N_2$, whence the curves $N_3,\ldots ,N_8$ are six independent
sections of the fibration. This fibration has two reducible fibers
of type $I_2$ (made up by the
classes $N_1$, $E_2-N_1$ and $N_2$, $E_2-N_2$).\\
The Segre map $s$ sends $\mathbb{P}^1\times\mathbb{P}^2$ in
$\mathbb{P}^5$.
$$
\begin{array}{ccccccc}
 &            &        &\mathbb{P}^1&&&\\
 &\phi_{L_2}&\nearrow&            &\nwarrow&&\\
X&            &        &\stackrel{\phi_{L_2}\times{\phi_{L_1}}}{\lra}&&\mathbb{P}^1\times\mathbb{P}^2&\stackrel{s}{\longrightarrow}\mathbb{P}^5\\
 &\phi_{L_1}&\searrow&&\swarrow&&\\
 &&&\mathbb{P}^2&&&
\end{array}
$$
Observe that the map $s\circ(\phi_{L_2}\times\phi_{L_1}):X\lra\mathbb{P}^5$ is
the map $\phi_{L_1+L_2}=\phi_{L-\hat{N}}$ (since $L_1+L_2=L-\hat{N}$). Indeed
let $s_1, s_2, s_3$ be a basis of $H^0(L_1)$, and $s_4, s_5$ be a basis of
$H^0(L_2)$. Then the products $s_1s_4$, $s_1s_5$, $s_2s_4$, $s_2s_5$, $s_3s_4$,
$s_3s_5$ are linear independent sections in $H^0(L_1+L_2)$ and define the Segre
embedding of $X$ in $\mathbb{P}^1\times\mathbb{P}^2$. Since $h^0(L_1+L_2)=6$
then the map $\phi_{L_1+L_2}$ is exactly the map
$s\circ(\phi_{L_2}\times\phi_{L_1})$.
\end{itemize}

\subsection{The case of $L^2=16$, $NS(X)=\cl_{16}'$, the map $\phi_{L_1}\times\phi_{L_2}$}\label{cl16'}

The intersection properties of $L_1$ and $L_2$ are $L_1\cdot L_1=2$, $L_2\cdot
L_2=2$ and $L_1\cdot L_2=4$. In \cite[Paragraph 3.6]{bertio} is described the
K3 surface $Y$ admitting a Nikulin involution with N\'eron Severi group
$\mathcal{M}_{8}$. Its quotient $\bar{X}$ is the complete intersection of a
hypersurface of bidegree $(1,1)$ and of a hypersurface of bidegree $(2,2)$ in
$\mathbb{P}^2\times \mathbb{P}^2$, its minimal resolution is $X$. A K3 surface
which is a complete intersection of a bidegree $(1,1)$ and a bidegree $(2,2)$
hypersurface in $\mathbb{P}^2\times\mathbb{P}^2$ is a {\it Wehler surface} (cf.
\cite{Wehler}). We describe this surface more in details in the Section
\ref{geometric}.
 The surface $X$ has a projective model in $\mathbb{P}^2\times
\mathbb{P}^2$ and the map associated to the divisors $L_1$ and
$L_2$ are respectively the projection to the first and to the
second projective space. The map
$\phi_{L_1}:X\rightarrow\mathbb{P}^2$ is a $2:1$ map. It contracts
the four rational curves $N_5,\ldots ,N_8$ to nodes of the branch
sextic of the double cover. The four curves $N_1,\ldots,N_4$ are
mapped to lines in $\mathbb{P}^2$. Since their intersection with
$L_1$ is equal to one, each of them is one of the two components
of the pullback of a line in $\mathbb{P}^2$. So their image under
the map $\phi_{L_1}$ is a line tritangent to the branch curve. The
branch curve has degree six and has four nodes so its genus is
$(6-1)(6-2)/2-4=6$. The curve $R_1$ on $X$ such that
$\phi_{L_1}(R_1)$ is the branch curve, has degree six and genus
six (because it is the branch curve, so the genus of the curve on
$X$ is the genus of its image on $\mathbb{P}^2$), this implies
that the curve $R_1$ has self-intersection ten ($g=R_1^2/2+1$) and
its intersection with $L_1$ is six. The curve $R_1$ has to
intersect the curves $N_i$, $i=1,\ldots, 4$ in three points
(because $N_i$ are mapped to tritangent to the sextic) and the
curves $N_i$, $i=5,\ldots, 8$ in two points (because the branching
curve has nodes in the points which are the images of these
curves). So we find
$$
R_1=3\left(\frac{L-N_1-N_2-N_3-N_4}{2}\right)-(N_5+N_6+N_7+N_8).
$$
Exactly in the same way one sees that the branch curve of the second projection
is
$$
R_2=3\left(\frac{L-N_5-N_6-N_7-N_8}{2}\right)-(N_1+N_2+N_3+N_4).
$$
The equation of a generic K3 surface which is complete intersection of a
$(1,1)$ and a $(2,2)$ hypersurface in
$\mathbb{P}^2\times\mathbb{P}^2$ is given by the system
$$
\left\{
\begin{array}{l}
\sum_{i,j=0,1,2}q_{ij}(x_0:x_1:x_2)y_iy_j=0\\
\sum_{i=0,1,2}l_i(x_0:x_1:x_2)y_i=0
\end{array}
\right.
$$
where $q_{ij}$ and $l_i$, $i,j=0,1,2$ are homogeneous polynomial
of degree respectively two and one in the variables $x_j$, which
are the coordinates of the first copy of $\mathbb{P}^2$ and $y_j$
denote coordinates of the second
copy of $\mathbb{P}^2$.\\
For a generic point
$(\overline{x_0}:\overline{x_1}:\overline{x_2})$ of $\mathbb{P}^2$
the system has two solutions in $(y_0:y_1:y_2)$ and this gives the
2:1 map to $\mathbb{P}^2$. If the point
$(\overline{x_0}:\overline{x_1}:\overline{x_2})$ is such that
$\sum_{i=0,1,2}l_i(\overline{x_0}:\overline{x_1}:\overline{x_2})y_i=0$
for each $(y_0:y_1:y_2)$, then the fiber on it is the quadric
$\sum_{i,j=0,1,2}q(\overline{x_0}:\overline{x_1}:\overline{x_2})y_iy_j=0$.
Otherwise if
$\sum_{i,j=0,1,2}q(\overline{x_0}:\overline{x_1}:\overline{x_2})y_iy_j=0$
for each $(y_0:y_1:y_2)$ then the fiber on
$(\overline{x_0}:\overline{x_1}:\overline{x_2})$ is a line.\\
Since in our case each map to $\mathbb{P}^2$ contracts four
lines, then for each copy of $\mathbb{P}^2$ there are four points
in which $\sum_{i,j=0,1,2}q(x_0:x_1:x_2)y_iy_j=0$ are identically
satisfied. Up to a projective transformation one can suppose that
the four points with dimension one fiber are, on each $\mathbb{P}^2$,
$(1:0:0), (0:1:0), (0:0:1), (1:1:1)$. This implies that the
equation $\sum_{i,j=0,1,2}q(x_0:x_1:x_2)y_iy_j=0$ is of the form
$$
\begin{array}{c}
y_0y_1(x_0x_1+ax_0x_2-(a+1)x_1x_2)+y_0y_2(bx_0x_1+cx_0x_2-(b+c)x_1x_2)+\\
+y_1y_2(-(1+b)x_0x_1-(a+c)x_0x_2+(a+b+c+1)x_1x_2)=0.
\end{array}
$$
and so it depends on three projective parameters. The equation of
type $(1,1)$ are
$$
\begin{array}{l}
x_0y_0+dx_0y_1+ex_0y_2+fx_1y_0+gx_1y_1+hx_1y_2+lx_2y_0+mx_2y_1+nx_2y_2=0
\end{array}
$$
and so depends on eight parameters (we can not apply other projective
transformations because we have chosen the points on which there are lines as
fibers). So a Wehler K3 surface such that the projection $\phi_{L_1}$ to
$\mathbb{P}^2$ contracts four rational curves of the K3 surface and
$\phi_{L_2}$ contracts four other rational curves (disjoint from the previous
curves) depends exactly on eleven parameters.

\subsection{The case of $L^2=24$, $NS(X)=\cl_{24}'$, the map $\phi_{L_1}\times\phi_{L_2}$ }\label{cl24'}

The intersection properties of $L_1$ and $L_2$ are $L_1\cdot
L_1=4$, $L_2\cdot L_2=4$ and $L_1\cdot L_2=6$. Each of them
defines a map from $X$ to $\mathbb{P}^3$ (by Lemma
\ref{lemma:birational2}). Each map $\phi_{L_i}$, $i=1,2$ contracts
four rational curves and sends the other in four lines. The curve
$L_1$
is sent by $\phi_{L_2}$ to a curve of degree six, and viceversa.\\
In \cite[Paragraph 3.8]{bertio}  it is described the K3 surface $Y$ admitting a
Nikulin involution with N\'eron-Severi group $\mathcal{M}_{12}$. Its quotient
is $\bar{X}$ and it is a complete intersection of four varieties of bidegree
$(1,1)$ in $\mathbb{P}^3\times \mathbb{P}^3$, the minimal resolution is $X$
(cf. Table \ref{tab1}). The projections to the two copies of $\mathbb{P}^3$ are
$\phi_{L_1}$ and $\phi_{L_2}$.


\section {Geometric conditions to have an even set.}\label{geometric}

In this section we describe geometrical properties of K3
surfaces which imply the presence of an even set. These even
sets can be of eight nodes, of eight rational curves (lines or
conics) or of some nodes and some rational curves. The following
results are in a certain sense the converse of the results of the
previous section, where we supposed that a K3 surface admits an
even set and we described its geometry. To prove the existence of
an even set on $S$ we will prove that either the lattice
$\cl_{2d}$ or $\cl_{2d}'$ is embedded in  $NS(S)$ and that the
sublattice $N$ of $\cl_{2d}$ (or of $\cl_{2d}'$) is generated over
$\mathbb{Q}$ by $(-2)$-irreducible curves. Since
rank $\cl_{2d}=$rank $\cl_{2d}'=9$ and since the K3 surfaces with N\'eron
Severi group equal to $\cl_{2d}$ or $\cl_{2d}'$ have an even
set, then the number of moduli of the families of K3 surfaces that we describe
here is eleven.

\subsection{Double cover of a cone with an even set.}

Let $S$ be a K3 surface which is a double cover of a cone, then by
\cite[Proposition 5.7 case iii)]{saintdonat} the map from $S$ to the cone is
given by a class $L'$ in $NS(S)$ such that either
\begin{itemize}
    \item[a)]  $L'=2E'+\Gamma_0+\Gamma_1$ with $\Gamma_0\cdot E'=\Gamma_1\cdot
E'=1$ and $\Gamma_0\cdot \Gamma_1=0$ or
    \item[b)]
    $L'=2E'+2\Gamma_0+\ldots+2\Gamma_n+\Gamma_{n+1}+\Gamma_{n+2}$,
    with $E'\cdot \Gamma_0=\Gamma_i\cdot\Gamma_{i+1}=1$ $i=0,\ldots,
    n-1$, $\Gamma_n\cdot\Gamma_{n+1}=\Gamma_{n}\cdot
    \Gamma_{n+2}=1$ and the other intersections are equal to zero.
\end{itemize}
The $\Gamma_i$'s are irreducible $(-2)$-curves. If we are in the case a), then
we can give a sufficient condition for $S$ to have an even set of eight
disjoint rational curves.
\begin{prop}\label{prop: even set double cover cone}
Let $S$ be a K3 surface such that there exist a map
$\phi_{L'}:S\stackrel{2:1}{\rightarrow} Z$, where $Z$ is a cone
and $L'=2E'+\Gamma_0+\Gamma_1$ with $\Gamma_0\cdot
E'=\Gamma_1\cdot E'=1$ and $\Gamma_0\cdot \Gamma_1=0$. If the
branch locus of the double cover is the union of a conic and a
sextic meeting in six distinct points and not passing through the
vertex of the cone, then $S$ admits an even set of eight disjoint
rational curves.
\end{prop}
\bprf We prove that under the hypothesis the lattice $\cl_4'$ is embedded in
the N\'eron Severi lattice of $S$. In particular there exist eight disjoint
rational curves in $NS(S)$ generating on $\Q$ a copy of $N$ in the N\'eron
Severi lattice. This implies that
$S$ admits an even set made up by these eight disjoint rational curves.\\
By the hypothesis the classes $L'$, $E'$ and $\Gamma_0$ are linearly
independent and  are in
$NS(S)$. The map $\phi_{L'}$ is a $2:1$ map to the cone, which
contracts the two rational curves $\Gamma_0$ and $\Gamma_1$ to the vertex of the cone. The
(smooth) K3 surface $S$ is the double cover of the blow up of the
cone in the vertex and in the six singular points of the
ramification locus. On $S$ there are six rational curves
$\Gamma_i$, $i=2,\ldots, 7$ on the six singular points of the
ramification locus, and the two rational curves $\Gamma_0$ and
$\Gamma_1$ on the blow up of the vertex of the cone (since this
 is not in the ramification locus we obtain two curves).\\
Let $C_2'$ be the curve such that $\phi_{L'}(C_2')=c_2'$ is the conic of the
branching locus. Since $c_2'$ is a conic, $C_2'\cdot L=2$ and since it does not
pass through the vertex then $C_2'\cdot\Gamma_0=C_2'\cdot\Gamma_1=0$, so
$C_2'\cdot E'=1$, moreover $C_2'$ is a rational curve and so $C_2'^2=-2$. Since
on the cone $c_2'$ passes through the six singular points of the ramification
locus, on the K3 surface we have $C_2'\cdot \Gamma_i=1$, $i=2,\ldots,7$.\\ The
classes $L'$, $E'$, $\Gamma_0$, $C_2$, $\Gamma_i$, $i=2,\ldots, 7$ spans a
lattice $R$ which is isometric to the lattice $\cl_4'$. In fact a basis for
$\cl_4'$ is given by $(L-N_1-N_2)/2$, $\hat{N}$ and $N_i$, $i=1,\ldots 7$. The
map
$$E'\mapsto (L-N_1-N_2)/2,\ \ \ C_2+E'-L'\mapsto \hat{N},\ \ \ \Gamma_i\mapsto
N_{i+1}\ \ i=0,\ldots, 7.$$ gives the explicit change of basis from $R$ to
$\cl_4'$. \eprf

\subsection{Complete intersection of one $(2,0)$ and three $(1,1)$
hypersurfaces in
$\mathbb{P}^4\times\mathbb{P}^2$}\label{subsection: ci P4*P2}

If $S$ is a complete intersection of a hypersurface of bidegree
$(2,0)$ and three hypersurfaces of bidegree $(1,1)$ in
$\mathbb{P}^4\times\mathbb{P}^2$, by the adjunction formula $S$ is
a K3 surface. The N\'eron Severi group of a generic K3 surface
which is a complete intersection of this type is generated by the
two divisors $D_1$ and $D_2$ associated to the two projections.
The family of the K3 surfaces of this type has Picard number two
and so it has 18 moduli. To give the complete description of the
N\'eron Severi group we compute the intersection $D_1\cdot D_2$.
We describe here how to find $D_1\cdot D_2$ as explained in
\cite[Section 5]{bert Brauer}. On the K3 surface the divisors
$D_1$ and $D_2$ correspond to the restriction to $S$ of the pull
back of the hyperplane section of $\mathbb{P}^4$, respectively of
$\mathbb{P}^2$. We put $h=\mathbb{P}^3\times \mathbb{P}^2$ and
$k=\mathbb{P}^4\times\mathbb{P}^1$. It is clear that
$h^3=\{point\}\times\mathbb{P}^2$, and so $h^4=0$ because in
$\mathbb{P}^4$ it corresponds to the intersection between a point
and a space. In the same way one computes that
$k^2=\PP^4\times\{point\}$ (intersection of two lines in
$\mathbb{P}^2$) and $k^3=0$, $h^3k^2=1$ ($\{point\}\times
\{point\}$). The hypersurface of bidegree $(2,0)$ corresponds to
$2h$ (has degree two with respect to the first factor, so with
respect to $h$, and zero with respect to the second factor, $k$)
and the hypersurfaces of bidegree $(1,1)$ correspond to the
divisor $h+k$. Since $X$ is the complete intersection of one
hypersurface of bidegree $(2,0)$ and three hypersurfaces of
bidegree $(1,1)$, $X$ corresponds in
$\mathbb{P}^4\times\mathbb{P}^2$ to the divisor $2h(h+k)^3$. We
want to compute $D_1\cdot D_2$ which is $h\cdot k$ restricted to
$2h(h+k)^3$. Then $D_1\cdot D_2$ is equal to $hk(2h)(h+k)^3$ in
the six dimensional space $\mathbb{P}^4\times \mathbb{P}^2$. The
terms $h^ik^j$ with $i+j=6$ correspond to the intersections of
codimension six and so are a finite number of points. The sum of
the coefficients of these
terms is exactly the number of points, so $D_1\cdot D_2=6$.\\
Hence the general K3 surface which is complete intersection of a $(2,0)$
hypersurface and three $(1,1)$ hypersurfaces in
$\mathbb{P}^4\times\mathbb{P}^2$ has N\'eron-Severi lattice isometric to
$$\begin{array}{c}
\{\mathbb{Z}^2,
\left[\begin{array}{cc}6&6\\6&2\end{array}\right]\}.
\end{array}
$$
This is a sublattice of the N\'eron Severi lattice of any $K3$ surface which is
a  complete intersection of a $(2,0)$ hypersurface and three $(1,1)$
hypersurfaces in $\mathbb{P}^4\times\mathbb{P}^2$.\\

\begin{prop}\label{primai.c.}
Let $S$ be a complete intersection of one hypersurface of bidegree $(2,0)$ and
three hypersurfaces of bidegree $(1,1)$ in $\mathbb{P}^4\times \mathbb{P}^2$.
Let $\phi_{A_1}$ and $\phi_{A_2}$ be the projections to the first and to the
second factor associated to the pseudo ample class $A_1$, with $A_1^2=6$, and
to the pseudo ample class $A_2$, with $A_2^2=2$. If there exist eight curves
$R_i$, $i=1,\ldots, 8$ such that $\phi_{A_1}$ contracts all these curves to
eight nodes of the image and $\phi_{A_2}$ sends these curves in lines on
$\mathbb{P}^2$, then $R_i$, $i=1,\ldots, 8$ form an even set.\end{prop} \bprf
The idea of the proof is similar to the proof of Proposition \ref{prop: even
set double cover cone} and is based on the presence of certain divisors in
$NS(S)$. The divisors $A_j$, $j=1,2$, $R_i$, $i=1,\ldots, 8$ are contained in
the N\'eron Severi group. Nine of these classes are linearly independent. The
lattice generated by $A_1,$ $A_2,$ $R_1,$ $R_2,$ $R_3,$ $R_4,$ $R_5,$ $R_6,$
$R_7$ is embedded in $NS(S)$ and a computation shows that it is isometric to
the lattice $\cl_6$. Since the lattice $\cl_6$ contains an even set, also in
the N\'eron Severi group of $S$ there is an even set made up by
$R_1,\ldots,R_7,2A_2-2A_1+R_1+\ldots+R_7$.\eprf \textbf{Remark} Observe that
Proposition \ref{primai.c.} gives a sufficient condition for a K3 surface
complete intersection in $\mathbb{P}^4$ to have an even set of nodes (or of
eight rational curves in the minimal resolution).

\subsection{Complete intersection of three quadrics in $\mathbb{P}^5$ with an even set of nodes}

We give two different sufficient conditions for a K3 surface in
$\mathbb{P}^5$ to have an even set of nodes. These two
possibilities correspond to the fact that the N\'eron Severi group
of such a K3 surface, with Picard number nine, is equal either to
the lattice $\cl_8$ or $\cl_8'$.

\begin{prop} Let $S$ be a K3 surface admitting two maps
$\phi_{A_1}$, $\phi_{A_2}$ associated to the pseudo ample class
$A_1$ with $A_1^2=8$ and to the ample class $A_2$ with $A_2^2=4$.
If there exist eight curves $R_i$, $i=1,\ldots, 8$ such that
$\phi_{A_1}$ contracts all these curves to eight nodes and
$\phi_{A_2}$ sends these curves to lines on the quartic in
$\mathbb{P}^3$, then $R_i$, $i=1,\ldots, 8$ form an even
set.\end{prop} \bprf One can prove that $\mathcal{L}_{8}$ is
primitively embedded in $NS(S)$ as in the Propositions \ref{prop:
even set double cover cone} and \ref{primai.c.}.\eprf

\begin{prop}\label{prop: even set otto'} Let $X$ be a K3 surface
in $\PP^5$ having eight nodes. These nodes form an even set if $X$
is the complete intersection of a smooth quadric and two quadrics,
which are singular in two planes $H\cong \PP^2$ and $K\cong
\PP^2$, $H\cap K=\emptyset$ and four of the points are contained
in $H$ and the other four in $K$.\end{prop} \bprf Let
$h_0=h_1=h_2=0$ and $k_0=k_1=k_2=0$ be the  equations defining $H$
resp. $K$ in $\mathbb{P}^5$, then we can write the equations of
the two cones as $h_0h_1-h_2^2=0$ and $k_0k_1-k_2^2=0$. The
quadrics $h_ik_j=0$, $i=0,1$, $j=0,1$ meet the K3 surface in two
curves $C_i$, $C_j$ with  multiplicity two, which passes through
four singular points, resp. to the other four. So
$2(C_i+C_j)\in|2L-(N_1+\ldots+N_8)|$, which shows that
$N_1+\ldots+N_8$ form an even set.\eprf

\subsection{Double covers of $\mathbb{P}^2$}

Here we consider two different K3 surfaces with an even set which
admit maps $2:1$ to $\mathbb{P}^2$. The first one is a Wehler
surface, the second one is not. In the first case the curves of
the even sets are contracted to singular points of the branch
locus or are sent to lines  of $\mathbb{P}^2$ which are tritangent
to the ramification locus, in the second case they are contracted
or sent to conics. Other double covers of $\PP^2$ with an even set
are described in \cite{barth2}.
\subsubsection{The first case: complete intersections of bidegree $(1,1)$, $(2,2)$
in $\mathbb{P}^2\times\mathbb{P}^2$.} The complete intersections
of bidegree $(1,1)$ and $(2,2)$ in
$\mathbb{P}^2\times\mathbb{P}^2$ are the Wehler surfaces. The
projections to the two copies of $\mathbb{P}^2$ are 2:1 maps. It
is known (but can also be computed as in Section \ref{subsection:
ci P4*P2}) that the N\'eron Severi group of the generic member of
this family is the two dimensional lattice
$$\begin{array}{c}
\{\mathbb{Z}^2,
\left[\begin{array}{cc}2&4\\4&2\end{array}\right]\}.
\end{array}
$$
The number of moduli of the family of the Wehler K3 surfaces is
18.

\begin{prop} Let $S$ be a Wehler K3 surface such that the
first projection $\pi_1$ contracts four rational disjoint curves
$R_l$, $l=1,\ldots, 4$ on $S$ and the second projection $\pi_2$
contracts other four rational disjoint curves $R_l$, $l=5,\ldots,
8$. Moreover the map $\pi_1$ sends the curves contracted by
$\pi_2$  to lines on $\mathbb{P}^2$ and viceversa. Then the eight
rational curves $R_l$, $l=1,\ldots, 8$ form an even set on
$S$.\end{prop} \bprf One can prove that $\mathcal{L}_{16}'$ is
primitively embedded in $NS(S)$ as in the Propositions \ref{prop:
even set double cover cone} and \ref{primai.c.}, and that
$N\subset\cl_{16}'$ is generated over $\Q$ by the curves $R_i$.\eprf
\subsubsection{The second case.}
\begin{prop} Let $S$ be a K3 surface admitting two maps $2:1$ to
$\mathbb{P}^2$. If there exist eight curves $R_i$, $i=1,\ldots, 8$
such that the map on the first copy of $\mathbb{P}^2$ contracts
the curves $R_i$, $i=1,\ldots, 4$ and sends the others in four
conics and the map on the second copy of $\mathbb{P}^2$ contracts
the curves $R_i$ $i=5,\ldots,8$ and sends the others in conics,
then $R_i$, $i=1,\ldots, 8$ is an even set of eight disjoint
rational curves.
\end{prop}
\bprf One can prove that $\mathcal{L}_{10}$ is primitively
embedded in $NS(S)$ as in the Propositions \ref{prop: even set
double cover cone} and \ref{primai.c.}.\eprf

\subsection{A mixed even set}
\begin{prop} Let $S$ be a K3 surface admitting two maps to
$\mathbb{P}^3$. If there exist eight curves $R_i$, $i=1,\ldots, 8$
such that the map on the first copy of $\mathbb{P}^3$ contracts
the curves $R_i$, $i=1,\ldots, 4$ and sends the others in four
conics and the map on the second copy of $\mathbb{P}^3$ contracts
the curves $R_i$, $i=5,\ldots,8$ and sends the others in conics,
then $R_i$, $i=1,\ldots, 8$ is an even set of eight disjoint
rational curves.
\end{prop}
\bprf One can prove that $\mathcal{L}_{12}$ is primitively embedded in $NS(S)$
as in the Propositions \ref{prop: even set double cover cone} and
\ref{primai.c.}.\eprf
In this case we have on a quartic in $\PP^3$ a {\it mixed
even set}, in fact it consists of four nodes and of four conics.

\subsection{Surfaces of bidegree $(2,3)$ in
$\mathbb{P}^1\times\mathbb{P}^2$.}

A K3 surface in $\mathbb{P}^1\times \mathbb{P}^2$ has bidegree $(2,3)$ by the adjunction formula.
These K3 surfaces are studied in \cite[Paragraph 5.8]{bert
Brauer}. The family has 18 moduli, in fact the N\'eron Severi
group of such a K3 surface has to contain two classes $D_1$, $D_2$
giving the regular maps $\phi_{D_1}$ and $\phi_{D_2}$, which
correspond to the projections to $\mathbb{P}^1$ and to
$\mathbb{P}^2$. Since $D_1$ defines a map to $\mathbb{P}^1$, we
have $D_1^2=0$ and analogously $D_2^2=2$, we compute the
intersection $D_1\cdot D_2$ as in the Paragraph \ref{primai.c.}.
On the K3 surface the divisors $D_1$ and $D_2$ correspond to the
hyperplane sections of $\mathbb{P}^1$, respectively of
$\mathbb{P}^2$. We put $h=\{point\}\times \mathbb{P}^2$ and
$k=\mathbb{P}^1\times\mathbb{P}^1$. It is clear that $h^2=0$,
$k^3=0$ and $hk^2=1$. Since the K3 surface has bidegree $(2,3)$ it
corresponds to the divisor $2h+3k$. Then $D_1\cdot D_2$
corresponds to $hk(2h+3k)$ in the three dimensional space
$\mathbb{P}^1\times \mathbb{P}^2$, so $D_1\cdot D_2=3$. We obtain
that the general K3 surface which has bidegree $(2,3)$ in
$\mathbb{P}^1\times\mathbb{P}^2$ has N\'eron-Severi lattice
isometric to
$$\begin{array}{c}
\{\mathbb{Z}^2,
\left[\begin{array}{cc}0&3\\3&2\end{array}\right]\}.
\end{array}
$$
This is a sublattice of the N\'eron Severi lattice of all the $K3$ surfaces
which
have  bidegree $(2,3)$ in $\mathbb{P}^1\times\mathbb{P}^2$.\\

\begin{prop}\label{prop: even set on P1*P2} Let $S$ be a K3 surface of
bidegree $(2,3)$ in $\mathbb{P}^1\times\mathbb{P}^2$ and such that the
projection to the first space $p_1$ (which gives an elliptic fibration)
contracts two disjoint rational curves and the projection $p_2$ to the second
space contracts other six disjoint rational curves. If the curves contracted by
$p_1$ are sent to lines by $p_2$ and the curves contracted by $p_2$ are sent by
$p_1$ to two sections of the elliptic fibration, then the eight rational curves
on $S$ form an even set.\end{prop} \bprf One can prove that $\mathcal{L}_{12}'$
is primitively embedded in $NS(S)$ as in the Propositions \ref{prop: even set
double cover cone} and \ref{primai.c.}.\eprf

\subsection{Complete intersections in $\mathbb{P}^3\times\mathbb{P}^3$.}

The N\'eron Severi group of a complete intersection of four hypersurfaces of
bidegree $(1,1)$ in $\mathbb{P}^3\times \mathbb{P}^3$ is generated by the two
divisors $D_1$ and $D_2$ associated to the two projections. The divisors $D_1$
and $D_2$ have self intersection equal to four, computing as before the
intersection between $D_1$ and $D_2$ one finds $D_1\cdot D_2=6$,
so the N\'eron Severi lattice of the generic K3 surface which is a complete
intersection of four bidegree $(1,1)$ hypersurfaces in $\mathbb{P}^3\times
\mathbb{P}^3$ is
$$\begin{array}{c}
\{\mathbb{Z}^2,
\left[\begin{array}{cc}4&6\\6&4\end{array}\right]\}.
\end{array}
$$

\begin{prop} Let $S$ be a complete intersection of four bidegree $(1,1)$ hypersurfaces
in $\mathbb{P}^3\times\mathbb{P}^3$. Let $A_1$ and $A_2$ be two
pseudo ample divisors defining two maps to $\mathbb{P}^3$. If the
map $\phi_{A_1}$, respectively $\phi_{A_2}$, contracts four
rational curves $R_1$, $R_2$, $R_3$, $R_4$, respectively $R_5$, $R_6$,
$R_7$, $R_8$, and sends the others four rational curves in lines,
then $R_i$, $i=1,\ldots, 8$ is an even set on $X$.\end{prop}
\bprf
One can prove that $\mathcal{L}_{24}'$ is primitively embedded in
$NS(S)$ as in the Propositions \ref{prop: even set double cover
cone} and \ref{primai.c.}.\eprf

\addcontentsline{toc}{section}{  \hspace{0.5ex} References}

\end{document}